\documentclass[11pt,a4paper]{amsart}
\usepackage[T1]{fontenc}
\usepackage{lmodern}
\usepackage[margin=1in]{geometry}
\usepackage[expansion=false]{microtype}
\usepackage{enumitem}
\usepackage{booktabs,array}
\usepackage{amssymb}
\usepackage{xcolor}
\usepackage[hidelinks]{hyperref}
\numberwithin{equation}{section}
\newtheorem{theorem}{Theorem}[section]
\newtheorem{proposition}[theorem]{Proposition}
\newtheorem{lemma}[theorem]{Lemma}
\newtheorem{corollary}[theorem]{Corollary}
\theoremstyle{definition}
\newtheorem{definition}[theorem]{Definition}

\newcommand{\R}{\mathbb R}
\newcommand{\C}{\mathbb C}
\newcommand{\Wcal}{\mathcal W}

\newcommand{\Ncal}{\mathcal N}
\newcommand{\Hone}{\dot H^1}
\newcommand{\dd}{\,\mathrm d}
\DeclareMathOperator{\Impart}{Im}
\DeclareMathOperator{\Repart}{Re}
\DeclareMathOperator{\tr}{tr}

\DeclareMathOperator{\dist}{dist}
\setlist[enumerate]{leftmargin=*,itemsep=3pt}
\allowdisplaybreaks[2]
\hypersetup{
 pdftitle={
 Dynamics of subcritical threshold solutions for the three-dimensional energy-critical NLS},
 pdfauthor={Meng Chen and Han Wang}
}
\title[Nonradial threshold classification]{Dynamics of subcritical threshold solutions for the three-dimensional energy-critical NLS}
\author{Meng Chen}
\author{Han Wang}

\address[Meng Chen, Han Wang]{Department of Mathematics and Statistics, Beijing Institute of Technology\\
Key Laboratory of Algebraic Lie Theory and Analysis of Ministry of Education\\
Beijing 100081, China}
\email[Meng Chen]{mchen@bit.edu.cn}
\email[Han Wang]{wang\_han@bit.edu.cn}
\date{}
\begin{document}
\raggedbottom
\begin{abstract}
We extend the radial threshold classification of Duyckaerts and Merle on the lower kinetic branch to nonradial solutions of the three-dimensional focusing energy-critical nonlinear Schr\"odinger equation in \(\dot H^1(\mathbb R^3)\). At the ground state energy and with initial gradient norm strictly below that of the ground state, every solution is global and either scatters in both time directions or is a symmetry transform of the special solution \(W^-\), possibly after conjugate time reversal. The proof combines Su and Zhao's below-threshold scattering theorem with a coercivity estimate for the interaction form at threshold under a kinetic gap and a local analysis near the ground state orbit.
\end{abstract}
\maketitle
\tableofcontents

\section{Introduction}
\label{sec:introduction}

\subsection{The equation and the ground state}
In this paper, we consider the focusing energy-critical equation
\begin{equation}\label{eq:nls}
 i\partial_tu+\Delta u+|u|^4u=0, 
 \qquad u(0)=u_0\in\Hone(\R^3).
\end{equation}
The conserved energy is
\[
 E(f)=\frac{1}{2}\int_{\R^3}|\nabla f|^2\dd x
          -\frac{1}{6}\int_{\R^3}|f|^6\dd x.
\]
The scaling $u(t,x)\mapsto\lambda^{\frac{1}{2}}u(\lambda^2t,\lambda x)$
preserves both the equation and the homogeneous energy norm.
A solution scatters forward if it is defined on $[0,\infty)$
and there is $u_+\in\Hone$ such that
\[
 \|u(t)-e^{it\Delta}u_+\|_{\Hone}\to0
 \qquad(t\to\infty).
\]
Backward scattering is defined in the same way. The time-reversal
symmetry used below is $u(t,x)\mapsto\overline{u(-t,x)}$.

The positive stationary solution
\begin{equation}\label{eq:ground-state}
 W(x)=\left(1+\frac{|x|^2}{3}\right)^{-\frac{1}{2}},
 \qquad -\Delta W=W^5,
\end{equation}
is an optimizer for the sharp Sobolev inequality
\cite{Aubin,Talenti}. Put
\begin{equation}\label{eq:basic-notation}
 G=\|\nabla W\|_2^2=\|W\|_6^6,\qquad E(W)=\frac{G}{3}.
\end{equation}
Then for any \(f \in \Hone\), 
\begin{equation}\label{eq:sharp}
 \|f\|_{L^6} \le G^{-\frac{1}{3}} \|\nabla f\|_{L^2}.
\end{equation}

The ground state orbit and its distance are
\begin{align}
 \Wcal&=\{\mathcal T_{\theta,\lambda,a}W:
                  \theta\in\R,\ \lambda>0,\ a\in\R^3\},
 \label{eq:ground-orbit}\\
 (\mathcal T_{\theta,\lambda,a}f)(x)
   &=e^{i\theta}\lambda^{\frac{1}{2}}f(\lambda(x-a)),\qquad
 d_{\Wcal}(f)=\inf_{w\in\Wcal}\|f-w\|_{\Hone}.
 \label{eq:symmetry-distance}
\end{align}
The parameters $\lambda$ and $\ell=\lambda^{-1}$ denote frequency
and physical scale, respectively.

Kenig and Merle \cite{KM} proved scattering below both ground state
thresholds for radial data in dimensions three through five.
The nonradial result was established by
Killip and Visan \cite{KV} in dimensions at least five and by
Dodson \cite{Dodson4} in dimension four. Su and Zhao \cite{SZ}
proved the three-dimensional result. Their scattering
theorem and frequency-localized interaction estimates provide the
global inputs to the present work. The scattering statement is
\begin{equation}\label{eq:below-scattering}
 E(f)<E(W),\quad \|\nabla f\|_{L^2}<\|\nabla W\|_{L^2}
 \quad\Longrightarrow\quad
 \text{the solution through $f$ scatters in both directions}.
\end{equation}
At energy $E(W)$, scattering is no longer the only possible behavior.
Duyckaerts and Merle \cite{DM} constructed a global radial solution
$W^-$ satisfying
\begin{equation}\label{eq:Wminus-properties}
 E(W^-)=E(W),\quad \|\nabla W^-(t)\|_{L^2}<\|\nabla W\|_{L^2},\quad
 \|W^-(t)-W\|_{\Hone}\lesssim e^{-e_0t}\quad(t\ge0),
\end{equation}
which scatters as $t\to-\infty$. The constant $e_0>0$ is the
positive real eigenvalue of the Hamiltonian linearization at $W$.

\subsection{Main result}

Duyckaerts and Merle \cite[Theorem 2]{DM} proved the radial
classification in dimensions three through five. They anticipated
its nonradial extension once the corresponding below-threshold
scattering theory became available \cite[Remark 1.4]{DM}. 
Building on the three-dimensional nonradial theory of Su and Zhao
\cite{SZ}, we establish this extension on the lower kinetic branch
in dimension three.

\begin{theorem}[Nonradial threshold classification]
\label{thm:classification}
Let $u$ be the maximal-lifespan strong solution of \eqref{eq:nls}
with
\[
 E(u_0)=E(W),\qquad\|\nabla u_0\|_{L^2} \le \|\nabla W\|_{L^2}.
\]
Then $u$ is global in both time directions. If $\|\nabla u_0\|_{L^2}  = \|\nabla W\|_{L^2}$, then
$u$ is a stationary element of $\Wcal$. If $\|\nabla u_0\|_{L^2}<\|\nabla W\|_{L^2}$, exactly
one of the following holds:
\begin{enumerate}
\item $u$ scatters in both directions;
\item $u$ is a fixed symmetry transform of $W^-$;
\item $u$ is a fixed symmetry transform of $\overline{W^-(-t)}$.
\end{enumerate}
Here a fixed symmetry transform of a solution $v$ means
\[
 u(t,x)=e^{i\theta_0}\lambda_0^{\frac{1}{2}}
 v\bigl(\lambda_0^2(t-t_0),\lambda_0(x-a_0)\bigr),
\]
with constant parameters. In particular, every solution in the strict
case scatters in at least one direction. No radial symmetry or
finite-mass assumption is imposed.
\end{theorem}

The equality case follows from \eqref{eq:sharp}; the strict lower
branch is treated in Section~\ref{sec:formal-classification}.

Table~\ref{tab:threshold-classification} summarizes the classification
results across dimensions. The case of equal gradient
norm follows from sharp Sobolev equality in every dimension.
For results on the upper kinetic branch, see \cite{Upper}.

\begin{table}[htbp]
\centering
\small
\setlength{\tabcolsep}{5pt}
\renewcommand{\arraystretch}{1.2}
\caption{Threshold classification in $\dot H^1(\mathbb R^d)$ under
$E(u_0)=E(W)$ and $\|\nabla u_0\|_2<\|\nabla W\|_2$,
where $W$ is the ground state in dimension $d$.
No finite-mass assumption is imposed.}
\label{tab:threshold-classification}
\begin{tabular}{@{}c
  >{\raggedright\arraybackslash}p{0.36\textwidth}
  >{\raggedright\arraybackslash}p{0.47\textwidth}@{}}
\toprule
Dimension & Radial data & General data\\
\midrule
$d=3$
& Duyckaerts--Merle \cite{DM}
& Present paper, Theorem~\ref{thm:classification}\\
$d=4$
& Duyckaerts--Merle \cite{DM}
& Ma--Miao--Murphy--Zheng \cite{MMMZ4}\\
$d\ge5$
& Duyckaerts--Merle \cite{DM} ($d=5$);
  Li--Zhang \cite{LZ} ($d\ge6$)
& Su--Zhao \cite{SZthreshold}\\
\bottomrule
\end{tabular}
\end{table}

The three-dimensional problem is naturally posed in $\Hone$:
$W(x)\sim |x|^{-1}$ at infinity, so even the ground state has
infinite mass. Finite mass and momentum are therefore unavailable
for general data in Theorem~\ref{thm:classification}. In the nonradial
setting, one must also allow a time-dependent spatial center.
These features make it necessary to control concentration without
assuming either a bounded center or a global $L^2$ bound. The
frequency-localized estimates of \cite{SZ} provide that control;
the compactness constructions of \cite{KV} supply the frequency
normalization needed to apply them.

\subsection{Outline of the proof}

We follow the compactness and local-dynamics framework of
Duyckaerts and Merle \cite{DM}, using the three-dimensional
nonradial estimates of Su and Zhao \cite{SZ}.
Below-threshold scattering and profile decomposition give compactness
of a nonscattering threshold orbit modulo scaling and translation.
At the threshold, however, the strict inequality $\|\nabla u(t)\|_{L^2}<\|\nabla W\|_{L^2}$ need
not have a uniform margin, as the solution $W^-$ shows. Thus one
cannot simply apply the below-threshold rigidity argument to the
whole orbit. We first exclude compact solutions that retain a fixed
kinetic gap; this yields approach to $\Wcal$ along a sequence of
times. The local dynamics then upgrade this sequential conclusion
to eventual trapping and classification.

For the global rigidity step, a compact tail with
$\|\nabla u(t)\|_{L^2}^2<\|\nabla W\|_{L^2}^2-\kappa$ produces a global compact model with the same
gap and a frequency bounded below. Even with this gap, the strict
directional coercivity used below threshold can fail at $E(W)$;
Section~\ref{subsec:kinetic-coercivity} gives an explicit example.
We instead use the cone of matrices generated by orthogonal
projections onto two-dimensional planes. The interaction Hessians
and their convolutions lie in this cone, and a two-dimensional
anisotropic dilation gives the required positive margin. This
margin persists in a neighborhood of the normalized orbit, allowing
for the small increase in energy that frequency truncation may cause.
The frequency-localized estimates then yield an integrable forcing
bound. The forward no-waste formula converts this bound into an
$L^2$ estimate tending to zero as $t\to\infty$, which contradicts
mass conservation.

Near $\Wcal$, modulation includes phase, scale and all three
translations. An integrable bound for the distance in rescaled time
controls the variation of these parameters. We combine hyperbolic
ejection, as in \cite{NR}, with an exponentially weighted comparison
against the special solutions constructed in \cite{DM} to identify
limits of normalized exits. Such a limit scatters in the outgoing
time direction, so half-line stability rules out exits after
arbitrarily close approaches. A nonscattering threshold solution is
therefore eventually trapped. The scalar modulation equations then
give exponential convergence to a fixed ground state, and the same
comparison identifies the solution with $W^-$.

Section~\ref{sec:preliminaries} collects the preliminaries.
Section~\ref{sec:threshold-kinetic-gap} proves rigidity under a
kinetic gap, and Section~\ref{sec:threshold-compactness} establishes
compactness and sequential approach to the ground state orbit.
Sections~\ref{sec:formal-local} and~\ref{sec:formal-classification}
develop the local dynamics and prove Theorem~\ref{thm:classification}.

\section{Preliminaries}
\label{sec:preliminaries}

\subsection{Notation and local theory}

Unless otherwise specified, all spatial norms and Sobolev spaces are over $\R^3$.
The notation $X\lesssim Y$ means $X\le CY$ for a fixed constant;
$X\sim Y$ means both $X\lesssim Y$ and $Y\lesssim X$.
Dependence on a fixed parameter is displayed when needed.
The homogeneous energy inner product is real:
\[
 \langle f,g\rangle_{\dot H^1}
   =\Repart\int\nabla f\cdot\overline{\nabla g}.
\]
We use the Fourier convention
$\widehat f(\xi)=\int e^{-2\pi i x\cdot\xi}f(x)\dd x$.
For a smooth radial cutoff $0\le\chi\le1$, equal to one on
$B(0,1)$ and supported in $B(0,2)$, put
\[
 \widehat{P_{\le\nu}f}(\xi)=\chi\bigl(\frac{\xi}{\nu}\bigr)\widehat f(\xi),
 \qquad P_{>\nu}=I-P_{\le\nu}.
\]
Dyadic projections are denoted by $P_N$. In particular the
high-frequency multiplier takes values in $[0,1]$.

For an interval $J$, write
\[
 S_J(u)=\int_{J\times\R^3}|u(t,x)|^{10}\dd x\dd t.
\]
A strong solution is a function $u\in C(I;\Hone)$ with
$S_J(u)<\infty$ for every compact $J\subset I$, satisfying the
Duhamel formula
\[
 u(t)=e^{i(t-t_0)\Delta}u(t_0)
       +i\int_{t_0}^t e^{i(t-s)\Delta}|u(s)|^4u(s)\dd s.
\]
We use the standard energy-space local theory developed in
\cite{CW,KM}. Solutions with initial data in $\Hone$ are unique on
their maximal lifespans and conserve the energy. A finite maximal
lifespan endpoint forces the scattering size to diverge, while a
global solution with finite scattering size scatters in the
corresponding time direction; see
\cite[Lemma~2.11 and Remark~2.15]{KM}. For initial data in
$H^1(\R^3)$, the usual persistence of regularity and conservation
of mass also hold.

For later use, define
\begin{equation}\label{eq:prelim-X}
 \|v\|_{X(J)}
 :=
 \|v\|_{L_t^\infty(J;\Hone)}
 +
 \|\nabla v\|_{L_t^2L_x^6(J)}.
\end{equation}
Interpolation and Sobolev embedding give
\begin{equation}\label{eq:X-embedding}
 \|v\|_{L^{10}_{t,x}(J)}
 \lesssim \|v\|_{X(J)}.
\end{equation}

We also use the long-time perturbation theorem
\cite[Theorem~2.14 and Remark~2.17]{KM}. In particular, if a
reference solution has bounded $L_t^\infty\Hone$ norm and finite
scattering size on an interval $J$, then solutions with sufficiently
close initial data exist throughout $J$ and remain close in $X(J)$.
For an approximate solution with equation error $e$, the corresponding
estimate also includes $\|\nabla e\|_{L_t^2L_x^{6/5}(J)}.$ 
These statements hold on bounded intervals and half-lines, in either
time direction, and require no finite-mass assumption. We shall use
them below for nonlinear profile approximations.

For $f\in C_c^\infty(\R^3)$, the dispersive
estimate gives
\begin{equation}\label{eq:free-L6-decay}
    \|e^{it\Delta}f\|_{L^6}
\le C|t|^{-1}\|f\|_{L^{6/5}}\to0
\quad (|t|\to\infty).
\end{equation}

The result extends to every $f\in\Hone$
by density and the uniform bound
\[
\|e^{it\Delta}f\|_{L^6}
\le C\|f\|_{\Hone}.
\]
Thus, if $u$ scatters in $\Hone$ as
$t\to\pm\infty$, then
\begin{equation}\label{eq:scatter-kinetic-limit}
\|u(t)\|_{L^6}^6\to0,\qquad
\|\nabla u(t)\|_{L^2}^2\to2E(u)
\end{equation}
in the corresponding time direction.

We also record some consequences of the sharp Sobolev inequality. For every $f\in\Hone$,
\begin{equation}\label{eq:positive-energy-lower}
\|\nabla f\|_{L^2}\le\|\nabla W\|_{L^2}
\quad\Longrightarrow\quad
E(f)\ge\frac13\|\nabla f\|_{L^2}^2.
\end{equation}
If $E(f)=E(W)$ and $\|\nabla f\|_{L^2}^2=G$, then \(\|f\|_{L^6}^6=G\).
Equality in \eqref{eq:sharp} therefore implies $f\in\Wcal$, including a constant complex phase.
Consequently, energy conservation, continuity, and uniqueness show that if
\[
E(u_0)\le E(W),\quad
\|\nabla u_0\|_{L^2}<\|\nabla W\|_{L^2},
\]
then \(\|\nabla u(t)\|_{L^2}<\|\nabla W\|_{L^2}\)
throughout the maximal lifespan. 
Indeed, if $\|\nabla u(t)\|_{L^2}
=\|\nabla W\|_{L^2}$ at some time,
energy conservation and the sharp Sobolev
inequality imply $u(t)\in\Wcal$.
By uniqueness, this contradicts the
strict inequality at the initial time.

\subsection{Profiles and compact normalized orbits}

We use the linear profile decomposition in $\Hone$
due to Keraani \cite{Keraani}, together with the
nonlinear profile approximation in \cite{KM,KV}.
Let $(f_n)$ be bounded in $\Hone$.
After passing to a subsequence, we have
\begin{equation}\label{eq:prelim-profile}
f_n=\sum_{j=1}^J
g_{\ell_n^j,x_n^j}e^{i\sigma_n^j\Delta}\phi^j
+w_n^J, \qquad (g_{\ell,x}f)(y)
=\ell^{-1/2}f\left(\frac{y-x}{\ell}\right).
\end{equation}
The scale, space, and time parameters are pairwise
asymptotically orthogonal, and each $\sigma_n^j$
is either zero or tends to $\pm\infty$.

Let $U^j$ denote the nonlinear profile associated
with $\phi^j$. If $\sigma_n^j=0$, then
$U^j(0)=\phi^j$. If $\sigma_n^j\to\pm\infty$,
then $U^j$ scatters to $e^{it\Delta}\phi^j$
as $t\to\pm\infty$, respectively.
In the latter case, \(E(U^j)=\frac{1}{2}\|\nabla\phi^j\|_{L^2}^2\)
by \eqref{eq:free-L6-decay}.

For each fixed $J$, we have
\begin{equation}\label{eq:prelim-gradient-decoupling}
\|\nabla f_n\|_{L^2}^2
=
\sum_{j=1}^J\|\nabla\phi^j\|_{L^2}^2
+\|\nabla w_n^J\|_{L^2}^2+o_n(1).
\end{equation}
\begin{equation}\label{eq:prelim-energy-decoupling}
E(f_n)
=
\sum_{j=1}^J E(U^j)
+E(w_n^J)+o_n(1),
\end{equation}
\begin{equation}\label{eq:prelim-small-remainder}
\lim_{J\to\infty}\limsup_{n\to\infty}
\|e^{it\Delta}w_n^J\|_{
L^{10}_{t,x}(\R\times\R^3)}=0.
\end{equation}

We apply nonlinear profile approximation only
on intervals where the relevant nonlinear
profiles have uniformly bounded scattering size.

\begin{definition}[Almost periodicity]
A solution is almost periodic on an interval $I$
if there are $N(t)>0$ and $x(t)\in\R^3$ such that
\[
\mathcal K=
\left\{
N(t)^{-1/2}
u\left(t,x(t)+N(t)^{-1}\cdot\right):
t\in I
\right\}
\]
is precompact in $\dot H^1$.
\end{definition}

We refer to $N$ as a frequency scale.
A phase parameter need not be included in
the definition since the phase group is compact.

Let $u$ be a nonzero maximal-lifespan almost
periodic solution, and write
$I=(T_-,T_+)$ for its maximal
lifespan. Compactness of the normalized orbit
and local well-posedness give the following
local constancy property
(see \cite[Lemma~1.5]{KV3}).
There exist $\delta>0$ and $C\ge1$ such that,
for every $t_0\in I$,
\[
[t_0-\delta N(t_0)^{-2},
 t_0+\delta N(t_0)^{-2}]
\subset I,
\]
and
\begin{equation}\label{eq:prelim-local-constancy}
C^{-1}N(t_0)\le N(t)\le CN(t_0)
\end{equation}
whenever $|t-t_0|\le\delta N(t_0)^{-2}.$ 
The constants are independent of $t_0$.
The argument applies to both the focusing
and defocusing equations.

\begin{lemma}
\label{lem:kin-finite-ap}
Let $u:I\times\R^3\to\C$ be a nonzero
maximal-lifespan strong solution to \eqref{eq:nls}
in $\dot H^1(\R^3)$.
If $u$ is almost periodic on $I$,
then $I=\R$.
\end{lemma}
\begin{proof}
Suppose its forward endpoint \(T\) is finite. Write its frequency scale
as \(N(t)\), physical scale as \(\ell(t)=N(t)^{-1}\), and normalized
precompact family as \(\mathcal K\subset\Hone\).
Uniform local existence on the compact family gives a positive
normalized lifespan. If \(N(t_n)\) stayed bounded along a sequence
\(t_n\to T\), with \(t_n<T\), stability would extend the solution beyond \(T\).
Thus \(N(t)\to\infty\) and \(\ell(t)\to0\).

Choose \(0\le\psi\le1\) in \(C_c^\infty(\R^3)\), equal to one
on \(B_1\) and supported in \(B_2\), and put
\(\psi_R(x)=\psi(\frac{x}{R})\).
For any \(\eta>0\), approximate \(\mathcal K\)
in \(\Hone\) to accuracy \(\eta\) by finitely many compactly supported
smooth functions. Their \(L^2\) norms have a finite maximum \(C_\eta\).
Critical scaling and local Sobolev embedding give, uniformly in the
translation parameters,
\[
 \|\psi_Ru(t)\|_2\le \ell(t)C_\eta+CR\eta.
\]
First let \(t\to T\) within the lifespan, then \(\eta\to0\). Hence
\(\|\psi_Ru(t)\|_2\to0\) for every fixed \(R\).

The localized continuity equation and Cauchy--Schwarz imply
\begin{equation}\label{eq:kin-sqrt-local-mass}
 \left|\|\psi_Ru(t)\|_2-\|\psi_Ru(s)\|_2\right|
 \le \frac{C}{R}\int_s^t \|\nabla u(\tau)\|_{L^2} \dd\tau.
\end{equation}
Indeed, the derivative of the squared local mass is
\(4\int\psi_R\nabla\psi_R\cdot
\Impart(\overline u\,\nabla u)\).
Regularizing the square root proves the displayed
estimate even when local mass vanishes.
Almost periodicity bounds \(\|\nabla u\|_{L^2}\) uniformly. Sending \(t\to T\)
in \eqref{eq:kin-sqrt-local-mass} gives
\[
 \|\psi_Ru(s)\|_2\le \frac{C}{R}(T-s)\sup_{\tau<T}\|\nabla u(\tau)\|_{L^2}.
\]
Letting \(R\to\infty\) forces \(u(s)=0\), a contradiction. The other
endpoint is identical. No global \(L^2\) norm was assumed or used.
\end{proof}

\section{Rigidity under a kinetic gap}
\label{sec:threshold-kinetic-gap}

The aim of this section is to rule out nonzero almost periodic solutions at energy \(E(W)\) whose gradient norm remains uniformly below that of \(W\). We first prove rigidity under a lower bound on the frequency scale. A compactness argument then removes this restriction, yielding Corollary 3.10.

\subsection{Plane coercivity}
\label{subsec:kinetic-coercivity}

A kinetic gap at $E(W)$ need not give a positive margin in every
direction. For $0<\alpha<1$, set
$f_\alpha(x)=\alpha W(\alpha^2x_1,x_2,x_3)$. Radiality of $W$ gives
\[
 \|\nabla f_\alpha\|_{L^2}^2=\frac{2+\alpha^4}{3}G,\qquad
\|f_\alpha\|_{L^6}^6=\alpha^4G,\qquad
 \|\partial_1f_\alpha\|_2^2=\frac{1}{3} \|f_\alpha\|_{L^6}^6.
\]
Thus $E(f_\alpha)=E(W)$ and
$\|\nabla f_\alpha\|_{L^2}^2<G$, but the
coercivity margin vanishes in the $x_1$ direction.
The following lemma shows that a uniform positive
margin nevertheless holds for every two-dimensional
plane, even when the energy slightly exceeds $E(W)$.
This is the form needed for the interaction Hessians
below.

\begin{lemma}[Plane coercivity]
\label{lem:tgap-robust-strip}
Let $0<\kappa<G$,
and suppose $f\in\Hone$ satisfies
\begin{equation}\label{eq:tgap-strip}
\|\nabla f\|_{L^2}^2
\le G-\kappa,
\quad
E(f)\le E(W)
+\frac{\kappa^2}{12G}.
\end{equation}
Then, for every two-dimensional plane
$\Pi\subset\R^3$,
\begin{equation}\label{eq:tgap-plane}
\|P_\Pi\nabla|f|\|_{L^2}^2
-\frac23\|f\|_{L^6}^6
\ge
\frac{\kappa}{6G}
\|f\|_{L^6}^6,
\end{equation}
where $P_\Pi$ denotes the orthogonal
projection onto $\Pi$.
\end{lemma}

\begin{proof}
The assertion is immediate if $f=0$.
Let $v=|f|$. By the diamagnetic inequality,
\[
 \|\nabla v\|_{L^2}^2
 \le \|\nabla f\|_{L^2}^2
 \le G-\kappa,
 \quad E(v)\le E(f).
\]
Moreover, $\|P_\Pi\nabla v\|_{L^2}>0$.
Indeed, otherwise $v$ would be constant along
$\Pi$, which is impossible for a nonzero
function in $L^6(\R^3)$.

Write $x=(y,z)\in\Pi\times\Pi^\perp$ and define \(v_\lambda(y,z)=\lambda v(\lambda y,z)\), here \( \lambda>0\).
A change of variables gives
\begin{equation}\label{eq:tgap-plane-scaling}
 \|\nabla v_\lambda\|_{L^2}^2
 =\lambda^2\|P_\Pi\nabla v\|_{L^2}^2
   +\|P_{\Pi^\perp}\nabla v\|_{L^2}^2, \qquad
 \|v_\lambda\|_{L^6}^6
 =\lambda^4\|v\|_{L^6}^6.
\end{equation}
Choose
\[
 \lambda^2
 =1+\frac{
 G-\|\nabla v\|_{L^2}^2
 }{\|P_\Pi\nabla v\|_{L^2}^2}>1.
\]
Then \(\|\nabla v_\lambda\|_{L^2}^2
 =G\).
By the sharp Sobolev inequality,
\(\lambda^4\|v\|_{L^6}^6
 \le G\).
Since
$\|v\|_{L^6}^3<\sqrt{G}$,
we obtain
\begin{equation}\label{eq:tgap-raw-plane}
 \|P_\Pi\nabla v\|_{L^2}^2
 \ge
 \frac{
 \bigl(G
       -\|\nabla v\|_{L^2}^2\bigr)
 \|v\|_{L^6}^3
 }{
 \sqrt{G}-\|v\|_{L^6}^3
 }.
\end{equation}
Set
\[
 d=\sqrt{G}-\|v\|_{L^6}^3>0,
 \quad
 e=(E(v)-E(W))_+.
\]
We have
\[
 G
 -\|\nabla v\|_{L^2}^2
 =
 \frac13\left(
 G-\|v\|_{L^6}^6
 \right)
 -2(E(v)-E(W)).
\]
Substituting this into
\eqref{eq:tgap-raw-plane} yields
\begin{equation}\label{eq:tgap-plane-energy-error}
 \|P_\Pi\nabla v\|_{L^2}^2
       -\frac23\|v\|_{L^6}^6\ge
 \frac{\|v\|_{L^6}^3}{3d}
 \left(d^2-6e\right).
\end{equation}
By the sharp Sobolev inequality and
\eqref{eq:tgap-strip},
\[
\begin{aligned}
 \|v\|_{L^6}^3
 \le
 \frac{\|\nabla v\|_{L^2}^3}
      {G}\le
 \frac{
 (G-\kappa)^{3/2}
 }{G}\le
 \frac{G-\kappa}
      {\sqrt{G}}.
\end{aligned}
\]
Hence \(d\ge\frac{\kappa}{\sqrt{G}}\).
Moreover, since $E(v)\le E(f)$, \(6e
 \le\frac{\kappa^2}
          {2G}
 \le\frac{d^2}{2}\).
It follows from
\eqref{eq:tgap-plane-energy-error} that
\[
\begin{aligned}
 \|P_\Pi\nabla v\|_{L^2}^2
       -\frac23\|v\|_{L^6}^6\ge\frac{d}{6}\|v\|_{L^6}^3\ge
 \frac{\kappa}
      {6\sqrt{G}}
 \|v\|_{L^6}^3\ge
 \frac{\kappa}
      {6G}
 \|v\|_{L^6}^6.
\end{aligned}
\]
Here the last inequality follows from
$\|v\|_{L^6}^3\le\sqrt{G}$.
Since $v=|f|$, this proves
\eqref{eq:tgap-plane}.
\end{proof}

\begin{corollary}
\label{cor:tgap-matrix-cone}
Under \eqref{eq:tgap-strip}, let \(P\) be a real symmetric matrix with
\begin{equation}\label{eq:tgap-cone}
 P\ge0,\qquad 2P\le(\tr P)I.
\end{equation}
Then
\begin{equation}\label{eq:tgap-constant-matrix}
 \int\nabla|f|^TP\nabla|f|
       -\frac{\tr P}{3} \|f\|_{L^6}^6
 \ge\frac{\kappa}{12G} \|f\|_{L^6}^6  (\tr P) .
\end{equation}
\end{corollary}

\begin{proof}
Let $e_1,e_2,e_3$ be an orthonormal
eigenbasis of $P$, with eigenvalues
$\lambda_1,\lambda_2,\lambda_3$.
Set $\Pi_j = e_j^\perp$ and $c_j = \frac{\operatorname{tr}P}{2} - \lambda_j$. 
Then,
\[
 \sum_{j=1}^3c_j
 =\frac{\operatorname{tr}P}{2},
 \quad
 P=\sum_{j=1}^3c_jP_{\Pi_j}.
\]
Applying Lemma~\ref{lem:tgap-robust-strip} to each $\Pi_j$, multiplying by $c_j$ and summing over
$j=1,2,3$, we get

\[
\begin{aligned}
\int_{\R^3}
\nabla|f|^T P\nabla|f|\,dx
-\frac{\operatorname{tr}P}{3}\|f\|_{L^6}^6=
\sum_{j=1}^3 c_j
\left(
\|P_{\Pi_j}\nabla|f|\|_{L^2}^2
-\frac23\|f\|_{L^6}^6
\right)\ge
\frac{\kappa}{6G}
\|f\|_{L^6}^6
\sum_{j=1}^3c_j.
\end{aligned}
\]
This proves the conclusion.
\end{proof}

\subsection{Static interaction coercivity}
The matrix cone in Corollary 3.2 contains the Hessians of the interaction kernel used below. We use this fact to obtain interaction coercivity uniformly near a compact family in \(\dot H^1\).
Use the kernel of \cite[Section 4]{SZ},
\[
 k_R(x)=|x|-R+Re^{-\frac{|x|}{R}},\qquad H_R=D^2k_R.
\]
Its radial and tangential eigenvalues are
\[
 a_R(r)=R^{-1}e^{-\frac{r}{R}},\qquad
 b_R(r)=\frac{1-e^{-\frac{r}{R}}}{r},
\]
with their continuous values at the origin. The inequality
\(1+s\le e^s\) shows that \(0\le a_R\le b_R\). Therefore
\begin{equation}\label{eq:tgap-hessian-cone}
 H_R\ge0,\qquad 2H_R\le(\Delta k_R)I.
\end{equation}

We use the interaction notation of \cite[Sections 2.3 and 4.1]{SZ}.
For real \(v\in H^1\), set
\[
 \begin{split}
 Q_R(v)&=-\iint\Delta^2k_R(x-y)v(x)^2v(y)^2\dd x\dd y,\\
 K_R(v)&=\iint v(y)^2\nabla v(x)^TH_R(x-y)\nabla v(x)\dd x\dd y,\\
 \mathcal N_R(v)&=\iint\Delta k_R(x-y)v(x)^6v(y)^2\dd x\dd y.
 \end{split}
\]
The principal interaction form is
\[
 \mathcal B_R(v)=Q_R(v)+4K_R(v)-\frac{4}{3}\mathcal N_R(v).
\]
For complex $v\in H^1$, set $\rho_v=|v|^2$,
$p_v=\Impart(\overline v\nabla v)$, and $Q_R(v)=Q_R(|v|)$.
Define
\begin{align*}
 \mathcal B_R(v)={}&Q_R(v)
 +4\iint \partial_j\partial_k k_R(x-y)
 \bigl[\rho_v(y)\Repart
        (\partial_j\bar v\,\partial_kv)(x)
           -(p_v)_j(y)(p_v)_k(x)\bigr]\dd x\dd y\\
 &-\frac{4}{3}\iint\Delta k_R(x-y)|v(x)|^6\rho_v(y)\dd x\dd y.
\end{align*}

For $R>0$, define
\[
\psi_R(z)
=
\left(1+\frac{|z|}{R}\right)e^{-|z|/R}
\]
and
\[
S_R(h)
=
\frac{1}{32\pi R^3}
\iint
\psi_R(x-y)|h(x)|^2|h(y)|^2
\,dx\,dy.
\]
Since
\[
\partial_R k_R(z)=\psi_R(z)-1,
\]
this density average will determine the
time-dependent radius used below.
By the Fourier multiplier calculation
in \cite[Lemma~4.1]{SZ},
\[
8\pi S_R(h)
\le Q_R(h)
\le 8\pi\|h\|_4^4.
\]
Moreover, $S_R(h)$ is nonincreasing in $R$.

Repeated spatial indices are summed. This agrees with the real
form above.

The bilinear diamagnetic inequality of
\cite[Lemma 4.1, Step 3]{SZ} gives
\begin{equation}\label{eq:tgap-modulus}
 \mathcal B_R(h)\ge\mathcal B_R(|h|),\qquad
 Q_R(h)=Q_R(|h|).
\end{equation}
The reduction is valid without a momentum or finite-mass hypothesis on the underlying solution.

\begin{proposition}
\label{prop:tgap-static}
Let \(\mathcal C\subset\Hone\) be compact, with \(0\notin\mathcal C\),
and assume
\[
 \sup_{z\in\mathcal C}E(z)\le E(W),\qquad
 \sup_{z\in\mathcal C}\|\nabla z\|_{L^2}^2\le G-\kappa
\]
for some \(0<\kappa<G\). There exist \(\varepsilon_0,R_*>0\),
depending only on \(\mathcal C\) and \(\kappa\), such that every
\(h\in H^1(\R^3)\) with
\(\inf_{z\in\mathcal C}\|h-z\|_{\Hone}<\varepsilon_0\) satisfies
\begin{equation}\label{eq:tgap-static}
 \mathcal B_R(h)\ge Q_R(h),\qquad R\ge R_*.
\end{equation}
\end{proposition}

\begin{proof}
Set \( b_*=\inf_{z\in\mathcal C}\|z\|_{L^6}^6>0\).
By continuity and compactness, we can choose
$\varepsilon_0>0$ such that every $h$ in the
$\varepsilon_0$-neighborhood of $\mathcal C$
satisfies
\[
 \|\nabla h\|_{L^2}^2
 \le G-\frac{\kappa}{2},
 \quad
 E(h)\le E(W)
 +\frac{\kappa^2}{48G},
 \quad
 \|h\|_{L^6}^6\ge\frac{b_*}{2}.
\]
Let $v=|h|$ and define
\[
 M_R=H_R*v^2,\quad
 m_R=\tr M_R,\quad
 m_0=m_R(0)>0.
\]
Since $v^2\ge0$ and \eqref{eq:tgap-hessian-cone},
we have
\[
 M_R(x)\ge0,\qquad
 2M_R(x)\le m_R(x)I.
\]
Thus Corollary~\ref{cor:tgap-matrix-cone},
applied with $\kappa/2$ to $M_R(0)$, gives
\begin{equation}\label{eq:tgap-frozen-margin}
\begin{aligned}
\int_{\R^3}\nabla v^T M_R(0)\nabla v\,dx
 -\frac{m_0}{3}\|v\|_{L^6}^6 \ge
 \frac{\kappa}{24G}
 \|v\|_{L^6}^6m_0
 \ge\sigma m_0,
\end{aligned}
\end{equation}
where \(\sigma=\frac{\kappa b_*}
 {48G}>0\).
This estimate remains valid if
$\varepsilon_0$ is decreased.

We now use the freezing and exterior
absorption estimates from
\cite[Proposition 4.3, Steps 2--3]{SZ}.
Fix a small $\delta>0$.
By compactness and Hardy's inequality,
we can choose $C_0\ge1$ and decrease
$\varepsilon_0$ so that
\[
 \int_{|x|>C_0}
 \left(
 |\nabla v|^2+\frac{v^2}{|x|^2}
 \right)\,dx\le\delta.
\]
Choose $0\le\chi\le1$ such that
$\chi=0$ on $B_{C_0}$,
$\chi=1$ outside $B_{2C_0}$, and
$|\nabla\chi|\le C/C_0$.
Set
\[
 v_e=\chi v,\qquad
 K_e=\int_{|x|>C_0}
 \nabla v^T M_R\nabla v\,dx.
\]
By the kernel derivative estimate
\cite[(4.3)]{SZ}, the error from replacing
$M_R(x)$ by $M_R(0)$ on $B_{2C_0}$
is $o_R(1)m_0$, uniformly in $h$.
Combining this estimate with
\eqref{eq:tgap-frozen-margin}, we obtain,
for $\delta\le\sigma/4$ and sufficiently
large $R$,
\begin{equation}\label{eq:tgap-static-preabsorb}
 \mathcal B_R(v)
 \ge Q_R(v)+2\sigma m_0+4K_e
 -\frac43\int_{\R^3}m_Rv_e^6\,dx.
\end{equation}

Next, write \(v^2=v_e^2+(1-\chi^2)v^2\).
Since $(1-\chi^2)v^2\ge0$ is supported
in $B_{2C_0}$, the weighted Sobolev
inequality \cite[Proposition 4.2]{SZ}
gives, for $R\ge2C_0$,
\[
 \int_{\R^3}m_Rv_e^6\,dx
 \le C_{\rm tail}\delta^2K_e
 +C_{\rm tail}\delta^3m_0.
\]
These estimates use the full density
$v^2$, with constants independent of
$\|h\|_{L^2}$. In particular, no energy
estimate for a spatial truncation is needed.

Choose $\delta>0$ sufficiently small
that
\[
 \delta\le\frac{\sigma}{4},\qquad
 \frac43C_{\rm tail}\delta^2\le2,\qquad
 \frac43C_{\rm tail}\delta^3\le\sigma.
\]
We then fix $C_0$ and $\varepsilon_0$,
and choose $R_*\ge2C_0$ sufficiently
large that the freezing error is at most
$\sigma m_0/4$ for every $R\ge R_*$.
It follows from
\eqref{eq:tgap-static-preabsorb} that
\[
 \mathcal B_R(v)
 \ge Q_R(v)+2K_e+\sigma m_0
 \ge Q_R(v).
\]
The result for complex-valued $h$
now follows from \eqref{eq:tgap-modulus}.
\end{proof}

\subsection{Frequency-localized estimates}
\label{subsec:audited-forward-no-waste}
Under the static coercivity hypothesis, our goal is to prove
\[
\int_0^\infty
\bigl\||u(t)|^4u(t)\bigr\|_{L^1}^2\,dt<\infty.
\]
The next two lemmas explain why this bound forces an almost periodic solution to vanish.
\begin{lemma}[Forward no-waste formula]
\label{lem:audited-forward-no-waste}
Let $u$ be a forward-global strong solution
that is almost periodic on $[0,\infty)$
in $\Hone(\R^3)$. 
Then, for every $t\ge0$,
\begin{equation}\label{eq:no-waste-formula}
    e^{i(t-T)\Delta}u(T)
\rightharpoonup0
\quad\text{in }\dot H^1(\R^3)
\quad(T\to\infty).
\end{equation}
Consequently,
\begin{equation}
    \label{eq:weak-limit}
     u(t)=-i\lim_{T\to\infty}
       \int_t^T e^{i(t-s)\Delta}|u(s)|^4u(s)\dd s,
\end{equation}
where the limit is weak in \(\Hone\).
\end{lemma}

\begin{proof}
Fix $t\ge 0$ and let $T_n\to\infty$.
By almost periodicity, we may write
\(u(T_n)=g_{N_n^{-1},x_n}f_n\),
where, after passing to a subsequence,
\[
 f_n\to f
 \quad\text{strongly in }\dot H^1(\R^3).
\]
Unitarity gives
\[
\begin{aligned}
\bigl\|e^{i(t-T_n)\Delta}u(T_n)
-e^{i(t-T_n)\Delta}g_{N_n^{-1},x_n}f
\bigr\|_{\dot H^1}=\|f_n-f\|_{\dot H^1}\to0,
\end{aligned}
\]
so it suffices to show
\begin{equation}\label{eq:no-waste-reduced}
 e^{i(t-T_n)\Delta}g_{N_n^{-1},x_n}f
 \rightharpoonup0
 \quad\text{in }\dot H^1.
\end{equation}

Passing to a further subsequence, assume that
$N_n\to0$, $N_n\to\infty$, or
$N_n\to N_*\in(0,\infty)$.
By density and unitarity, it suffices to
test \eqref{eq:no-waste-reduced} against
$\varphi\in\mathcal S(\R^3)$, with both
$f$ and $\varphi$ having compact Fourier
support away from the origin.

If $N_n\to0$ or $N_n\to\infty$, the Fourier
supports of $g_{N_n^{-1},x_n}f$ and
$\varphi$ are disjoint for all large $n$.
Hence their $\dot H^1$ inner product
vanishes, even after free evolution.

If $N_n\to N_*\in(0,\infty)$, let
$h=g_{N_*^{-1},0}f$.
By continuity of scaling, \(\|g_{N_n^{-1},x_n}f
 -h(\cdot-x_n)\|_{\dot H^1}\to0\).
The dispersive estimate gives
\[
\begin{aligned}
\left|
\left\langle
e^{i(t-T_n)\Delta}h(\cdot-x_n),
\varphi
\right\rangle_{\dot H^1}
\right|\le C|T_n-t|^{-3/2}
\|\nabla h\|_{L^1}\|\nabla\varphi\|_{L^1}
\to0,
\end{aligned}
\]
uniformly in $x_n$.
This proves \eqref{eq:no-waste-reduced}.
Since every sequence $T_n\to\infty$
has a subsequence satisfying this
convergence, we obtain \eqref{eq:no-waste-formula}.
By the Duhamel formula, letting $T\to\infty$, we get \eqref{eq:weak-limit}.
\end{proof}

\begin{lemma}
\label{lem:source-no-waste-zero}
Let \(u\) be a forward-global strong solution of quintic NLS in
\(\R^3\), with either sign of the nonlinearity. Suppose
\[
 M:=\sup_{t\ge0}\|u(t)\|_{\Hone}<\infty,
 \qquad
 \int_0^\infty\bigl\||u(s)|^4u(s)\bigr\|_1^2\dd s<\infty,
\]
and assume the forward no-waste Duhamel formula in \(\Hone\).
Then \(u\equiv0\). 
\end{lemma}
\begin{proof}
Put \(F=|u|^4u\) and
\(Q_t=\int_t^\infty\|F(s)\|_1^2\dd s\). For each \(L>0\),
dual endpoint Strichartz \cite{KT} and Bernstein give
\[
 \Big\|\int_t^T e^{i(t-s)\Delta}P_{\le L}F(s)\dd s\Big\|_2
 \le C\|P_{\le L}F\|_{L_s^2L_x^{\frac{6}{5}}((t,T))}
 \le CL^{\frac{1}{2}}Q_t^{\frac{1}{2}}.
\]
The truncated integrals are Cauchy in \(L^2\) as \(T\to\infty\).
Their limit agrees distributionally with \(P_{\le L}u(t)\)
by no-waste: weak convergence in \(\Hone\) implies convergence
against Schwartz tests, using the embedding in \(L^6\).
Consequently
\[
 \|P_{\le L}u(t)\|_2\le CL^{\frac{1}{2}}Q_t^{\frac{1}{2}},
 \qquad \|P_{>L}u(t)\|_2\le C\frac{M}{L}.
\]
In particular \(u(t)\in L^2\) at every finite time. Optimizing
with \(L=M^{\frac{2}{3}}Q_t^{-\frac{1}{3}}\), when both quantities are positive,
yields
\begin{equation}\label{eq:source-tail-mass}
 \|u(t)\|_2\le CM^{\frac{1}{3}}Q_t^{\frac{1}{3}}\to0.
\end{equation}
The zero cases follow by taking a limit in the preceding two bounds.
Energy-space uniqueness and \(H^1\) persistence identify \(u\)
with the finite-mass solution through any finite time. Its mass
is therefore conserved, and \eqref{eq:source-tail-mass} forces it
to be zero.
\end{proof}

We record the frequency-localized estimates of \cite{SZ} separately
from the static coercivity established above. For a forward-global
solution $u$ with bounded gradient, write
\[
 F(z)=|z|^4z,\qquad (i\partial_t+\Delta)u=-F(u),\qquad
 B_0=\sup_{t\ge0}\|\nabla u(t)\|_2.
\]
For $\nu>0$, use the notation
\[
 \begin{gathered}
 l=P_{\le\nu}u,\quad h=P_{>\nu}u,\quad a(t)=\|h(t)\|_4^4,\\
 \eta_\nu=\sup_{t\ge0}\|\nabla l(t)\|_2,\quad
 M_\nu=\sup_{t\ge0}\|h(t)\|_2,\quad d_\nu=\nu M_\nu.
 \end{gathered}
\]
The projected equation is
$(i\partial_t+\Delta)h=-F(h)-\mathcal E$, where
\[
 g=P_{\le\nu}F(u),\qquad
 \mathcal E=P_{>\nu}F(u)-F(h),\qquad
 \rho=|h|^2,\qquad p=\Impart(\bar h\nabla h).
\]
With the mass and momentum brackets
\[
 \{z,h\}_m=\Impart(z\bar h),\qquad
 \{z,h\}_p=\Repart(z\nabla\bar h-h\nabla\bar z),
\]
the source contribution associated with a kernel $k$ is
\begin{equation}\label{eq:transfer-source-definition}
\begin{split}
 \operatorname{Err}_{k,\nu}(t)
 =-2\iint\nabla k(t,x-y)\cdot
 \bigl[\rho(t,y)\{\mathcal E,h\}_p(t,x)+2p(t,x)\{\mathcal E,h\}_m(t,y)\bigr]\dd x\dd y.
\end{split}
\end{equation}

\begin{lemma}
\label{lem:energy-free-analytic-inputs}
Let $u:[0,\infty)\times\R^3\to\C$ be a
forward-global strong solution to \eqref{eq:nls}.
Suppose that there exist $N(t)\ge1$ and
$x(t)\in\R^3$ such that
\[
 \mathcal C=
 \overline{\left\{
 N(t)^{-1/2}u\left(t,x(t)+N(t)^{-1}\cdot\right):
 t\ge0
 \right\}}
\]
is compact in \(\Hone\) and $0\notin\mathcal C$.
With the notation above, set
\[
\begin{aligned}
 f(t)=\|F(u(t))\|_{L^1},\quad
 F_+(t)=\sup_{s>0}\frac1s
 \int_t^{t+s}f(r)\dd r,
\end{aligned}
\]
and, for $T>0$,
\[
\begin{aligned}
 E_f(T)=\int_0^T f(t)^2\dd t,\quad
 H_f(T)=\int_0^T F_+(t)^2\dd t,\quad
 K_a(T)=\int_0^T a(t)\dd t.
\end{aligned}
\]
There exists $\nu_0>0$ such that,
for every $0<\nu\le\nu_0$, the following
estimates hold.

\begin{enumerate}
\item As $\nu\to0$, \(\eta_\nu\to0\), \(d_\nu\to0\).
Moreover, for some $c_{\mathcal C}>0$
independent of $\nu$,
\begin{equation}\label{eq:transfer-concentration}
 N(t)a(t)\ge c_{\mathcal C},
 \quad t\ge0.
\end{equation}
We also have
\[
 \sup_{t\ge0}\|u(t)\|_{L^4}<\infty,
 \quad
 \sup_{t\ge0}f(t)<\infty.
\]

\item For almost every $t\ge0$,
\begin{equation}\label{eq:transfer-low}
 \|l(t)\|_{L^\infty}\le C\nu F_+(t),
 \quad
 \|\nabla l(t)\|_{L^\infty}
 \le C\nu^2F_+(t).
\end{equation}
For every $T>0$,
\begin{equation}\label{eq:transfer-source}
 E_f(T)\le C_uK_a(T)
 +C_u\eta_\nu^8H_f(T).
\end{equation}

\item Let $k(t,z)$ be real, even in $z$,
and measurable in $t$, with weak spatial
derivatives satisfying
\begin{equation}\label{eq:transfer-kernel-class}
 |\nabla k(t,z)|
 +|z||D^2k(t,z)|
 +|z|^2|D^3k(t,z)|
 \le C_k
\end{equation}
almost everywhere, uniformly in $t$.
Then, for every $T>0$,
\begin{equation}\label{eq:transfer-error}
\begin{aligned}
 \left|\int_0^T
 \operatorname{Err}_{k,\nu}(t)\dd t\right|\le C_{u,C_k}d_\nu^{1/2}
 \left[\nu^{-3}+K_a(T)+H_f(T)\right].
\end{aligned}
\end{equation}
\end{enumerate}
All constants are independent of $\nu$
and $T$. In \textup{(iii)}, the constant
depends on $k$ only through $C_k$.
\end{lemma}

\begin{proof}
By compactness and scaling invariance,
\[
 \sup_{t\ge0}\|\nabla u(t)\|_{L^2}
 \le\sup_{z\in\mathcal C}
 \|\nabla z\|_{L^2}<\infty.
\]
The frequency estimates
$\eta_\nu,d_\nu\to0$ follow from the
compactness argument in
\cite[Proposition 2.3, proof of (2.6)]{SZ}.

Since $0\notin\mathcal C$, there exist
$R_0,c_0>0$ such that \(\inf_{z\in\mathcal C}
 \int_{|x|\le R_0}|z(x)|^4\dd x
 \ge c_0\).
The uniform low-frequency smallness
then gives \eqref{eq:transfer-concentration}
by the argument of
\cite[Proposition 2.3, proof of (2.7)]{SZ}.
These estimates use compactness and
$N(t)\ge1$, but no energy threshold.

By Lemma~\ref{lem:audited-forward-no-waste},
$u$ satisfies the forward Duhamel
representation. Hence
\cite[Lemma 3.1]{SZ} gives the uniform
$L^4$ and $L^1$ bounds in \textup{(i)}.
The estimates in \textup{(ii)} follow from
\cite[Proposition 3.2]{SZ}.

Finally, we adapt the proof of
\cite[Proposition~6.1]{SZ}.
The additional dispersive input is the
long-time Strichartz estimate
\[
\|\nabla P_{\le L}u\|_{L_t^2L_x^6(I)}
\lesssim_u
1+L^{3/2}
\left(\int_I N(t)^{-1}\dd t\right)^{1/2},
\]
which follows from
\cite[Theorem~4.1 and (4.5)]{KV3}
under the uniform $\dot H^1$ bound and
almost periodicity. Local constancy and
the uniform local Strichartz estimates
extend this bound to arbitrary compact
intervals $I\subset[0,\infty)$.

By \eqref{eq:transfer-concentration},
\[
\int_0^T N(t)^{-1}\dd t
\lesssim_{\mathcal C}K_a(T).
\]
Set
\[
h=P_{>\nu}u,\quad
g=P_{\le\nu}(|u|^4u),\quad
Z=\nabla\Delta^{-1}h.
\]
Dyadic summation then gives
\[
\begin{aligned}
\|Z\|_{L_t^2L_x^6([0,T])}
\lesssim_u
\nu^{-2}
+\nu^{-1/2}K_a(T)^{1/2}\lesssim_u
\nu^{-1/2}
\bigl(\nu^{-3}+K_a(T)\bigr)^{1/2}.
\end{aligned}
\]
Moreover, Bernstein's inequality yields
\[
\|D^mg\|_{L_t^2L_x^p([0,T])}
\lesssim
\nu^{m+3(1-1/p)}
\||u|^4u\|_{L_t^2L_x^1([0,T])}
\]
for $m=0,1,2$ and $1\le p\le\infty$.

The convolution-field estimates and
nonlinear cancellations in
\cite[Proposition~6.1]{SZ}
apply under the kernel assumptions
stated above. They use the uniform
$\dot H^1$ bound, Hardy and
Hardy--Littlewood--Sobolev inequalities,
\eqref{eq:transfer-low},
\eqref{eq:transfer-source},
and the estimates just established.

All integrations by parts involving the
kernel are spatial, so measurability
in time is sufficient for the error
estimate. No energy threshold or
minimality assumption is needed in
these bounds. The argument of
\cite[Proposition~6.1]{SZ} therefore
yields \eqref{eq:transfer-error},
with constants independent of
$\nu$, $T$, and the individual kernel
in the stated class.
\end{proof}

\begin{proposition}
\label{prop:static-to-source}
Let $u$ satisfy Lemma~\ref{lem:energy-free-analytic-inputs}.
Assume that there are $\varepsilon_0,R_*>0$ such that
\begin{equation}\label{eq:transfer-static-hyp}
 \mathcal B_R(v)\ge Q_R(v)
 \quad\text{whenever }v\in H^1,\quad
 \dist_{\Hone}(v,\mathcal C)<\varepsilon_0,
 \quad R\ge R_*.
\end{equation}
Then, for some fixed $\nu>0$,
\begin{equation}\label{eq:transfer-global-integrability}
 \int_0^\infty\|F(u(t))\|_1^2\dd t
 +\int_0^\infty\|P_{>\nu}u(t)\|_4^4\dd t<\infty.
\end{equation}
\end{proposition}

\begin{proof}
We follow \cite[Sections 5--7]{SZ}, using
\eqref{eq:transfer-static-hyp} in place of
the static coercivity estimate in
\cite[Proposition 4.3]{SZ}.
Set $h=P_{>\nu}u$ and $a(t)=\|h(t)\|_4^4$.
For $\Lambda\ge1$, let $R(t)>0$ satisfy
\begin{equation}\label{eq:transfer-radius}
R(t)=\Lambda S_{R(t)}(h(t)).
\end{equation}
By \cite[Lemmas 4.1 and 5.2]{SZ},
\[
\kappa(\Lambda)a(t)\le R(t)\le\Lambda a(t),
\quad
\kappa(\Lambda)\to\infty
\quad(\Lambda\to\infty).
\]
By \eqref{eq:transfer-concentration},
$N(t)a(t)\ge c_{\mathcal C}>0$ for all
sufficiently small $\nu$. Hence
\[
N(t)R(t)\ge
\kappa(\Lambda)c_{\mathcal C}.
\]
Fix $\Lambda$ sufficiently large that
$N(t)R(t)\ge R_*$ for all $t\ge0$.

Let
\[
\widetilde u(t)
=g_{N(t)^{-1},x(t)}^{-1}u(t),
\quad
\widetilde h(t)
=g_{N(t)^{-1},x(t)}^{-1}h(t).
\]
Since $\widetilde u(t)\in\mathcal C$,
Lemma~\ref{lem:energy-free-analytic-inputs}
\textup{(i)} gives
\[
\begin{aligned}
\dist_{\Hone}(\widetilde h(t),\mathcal C)
\le
\|\widetilde h(t)-\widetilde u(t)\|_{\Hone}=\|\nabla P_{\le\nu}u(t)\|_2
\le\eta_\nu\to0.
\end{aligned}
\]
Choose $\nu$ sufficiently small that
$\eta_\nu<\varepsilon_0$.
Applying \eqref{eq:transfer-static-hyp}
to $\widetilde h(t)$ at radius
$N(t)R(t)$, and then scaling back,
we obtain
\begin{equation}\label{eq:transfer-principal}
\begin{aligned}
\mathcal B_{R(t)}(h(t))
\ge Q_{R(t)}(h(t))\ge\frac{8\pi R(t)}{\Lambda}
\ge c_u a(t).
\end{aligned}
\end{equation}

To integrate \eqref{eq:transfer-principal},
we use the interaction Morawetz action
\[
\mathcal M_R(t)
=
2\iint_{\R^3\times\R^3}
\rho(t,y)\nabla k_R(x-y)\cdot p(t,x)
\dd x\dd y.
\]
For the radius $R(t)$ chosen above, set
\[
\mathcal M(t)=\mathcal M_{R(t)}(t),
\quad
D(t)=R'(t)
\left.\partial_R\mathcal M_R(t)
\right|_{R=R(t)},
\]
where the $R$ derivative is taken with
$h(t)$ fixed.

The frequency-localized interaction
Morawetz identity gives
\begin{equation}
\label{eq:transfer-morawetz}
\begin{aligned}
\mathcal M(T)-\mathcal M(0)
=\int_0^T\mathcal B_{R(t)}(h(t))
\dd t+\int_0^T
\operatorname{Err}_{k_{R(t)},\nu}(t)
\dd t+\int_0^T D(t)\dd t.
\end{aligned}
\end{equation}
The identity is justified at energy
regularity by the regularization
argument of \cite[Section 2.3]{SZ}.

Since $\|\nabla k_R\|_\infty\le1$,
Bernstein's inequality gives
\begin{equation}
\label{eq:transfer-action-bound}
\begin{aligned}
|\mathcal M(t)|
\le 2\|h(t)\|_2^3\|\nabla h(t)\|_2\le C_u\nu^{-3}.
\end{aligned}
\end{equation}

The moving-radius estimate
\cite[Lemma 5.4]{SZ} gives
\begin{equation}\label{eq:transfer-radius-error}
D(t)\ge-C_u(\eta_\nu^2+\nu)a(t).
\end{equation}
By decreasing $\nu$, this error is
absorbed into \eqref{eq:transfer-principal}.

By \eqref{eq:transfer-principal},
\eqref{eq:transfer-radius-error},
and \eqref{eq:transfer-morawetz},
\[
\begin{aligned}
c_uK_a(T)
\le 
C_u\nu^{-3}
+C_u(\eta_\nu^2+\nu)K_a(T)+\left|
\int_0^T
\operatorname{Err}_{k_{R(t)},\nu}(t)
\dd t
\right|.
\end{aligned}
\]
Decreasing $\nu$ and using
\eqref{eq:transfer-error}, we obtain
\[
K_a(T)
\le C_u\nu^{-3}
+C_ud_\nu^{1/2}
\bigl(K_a(T)+H_f(T)\bigr).
\]

Since $d_\nu\to0$, we absorb the term
involving $K_a(T)$. Combining the
result with \eqref{eq:transfer-source},
we obtain
\[
E_f(T)\le A_\nu+\varepsilon_\nu H_f(T),
\]
where $A_\nu=C_u\nu^{-3}$ and
$\varepsilon_\nu
=C_u(d_\nu^{1/2}+\eta_\nu^8)\to0$.

To close the estimate, set
$b(T)=\sup_{S\ge T}E_f(S)/S<\infty$.
The one-sided maximal-function argument
of \cite[Section 7]{SZ} gives
$H_f(T)\le C_{\rm max}Tb(T)$.
Consequently,
\[
b(T)\le\frac{A_\nu}{T}
+C_{\rm max}\varepsilon_\nu b(T).
\]
Choose $\nu$ sufficiently small that
$C_{\rm max}\varepsilon_\nu\le1/2$.
Then $E_f(T)\le Tb(T)\le2A_\nu$
uniformly in $T$.

The global $L^2$ maximal inequality
now gives $H_f(\infty)<\infty$.
Returning to the estimate for $K_a(T)$,
we also obtain $K_a(\infty)<\infty$.
This proves
\eqref{eq:transfer-global-integrability}.
\end{proof}

\begin{proposition}
\label{prop:tgap-forward-rigidity}
Let $u:[0,\infty)\times\R^3\to\C$ be a
forward-global almost periodic strong solution
to \eqref{eq:nls}, with a frequency scale
$N(t)\ge1$ for $t\ge0$.
Suppose that
\[
 E(u)\le E(W),\quad
 \sup_{t\ge0}\|\nabla u(t)\|_{L^2}^2
 \le G-\kappa
\]
for some $\kappa>0$. Then $u\equiv0$.
\end{proposition}

\begin{proof}
Suppose $u\not\equiv0$, and let $\mathcal C$ be
the compact closure of the normalized forward orbit.
By the sharp Sobolev inequality and the kinetic gap, \(E(u)\ge\frac13\|\nabla u_0\|_{L^2}^2>0\).
Thus $0\notin\mathcal C$, and scaling invariance
shows that $\mathcal C$ satisfies the assumptions
of Proposition~\ref{prop:tgap-static}.
Combining Proposition~\ref{prop:tgap-static},
Lemmas~\ref{lem:audited-forward-no-waste}
and \ref{lem:energy-free-analytic-inputs},
and Proposition~\ref{prop:static-to-source},
we obtain \(\int_0^\infty
\||u(t)|^4u(t)\|_{L^1}^2\dd t<\infty\).
Lemma~\ref{lem:source-no-waste-zero}
then gives $u\equiv0$, a contradiction.
\end{proof}

\subsection{Frequency normalization}
Proposition 3.8 assumes \(N(t)\ge1\). We remove this restriction by passing to a normalized time-translate limit whose frequency scale is bounded below.

Let $w$ be an almost periodic solution
with frequency scale $N_w$ and spatial
center $x_w$, as in
Section~\ref{sec:preliminaries}.
 Given $t_n\in I(w)$, its normalized time translates are
\[
 w_n(s,y)=N_w(t_n)^{-\frac{1}{2}}
 w\left(t_n+\frac{s}{N_w(t_n)^2},\,
        x_w(t_n)+\frac{y}{N_w(t_n)}\right).
\]
A normalized time-translate limit is a maximal-lifespan strong
solution $u$, with $0\in I(u)$, such that, along a subsequence,
\[
 w_n\to u\quad\text{in }C(J;\Hone)
 \quad\text{for every compact }J \Subset I(u),
\]
with $J\subset I(w_n)$ for all sufficiently large $n$.
For the conjugate time reversal $w(t)=\overline{v(-t)}$, we use
$N_w(t)=N_v(-t)$ and $x_w(t)=x_v(-t)$.

\begin{lemma}
\label{lem:closed-class-floor}
Let $v$ be a nonzero global almost periodic
strong solution. Assume that every normalized
time-translate limit of $v$ or
$\overline{v(-t)}$ is nonzero, global, and
almost periodic with a common compactness modulus.
Then, up to a fixed critical scaling,
one such limit $u$ has a frequency function
$N_u$ satisfying \(N_u(t)\ge1\), \(t\in\R\).
\end{lemma}

\begin{proof}
We follow the frequency-normalization argument
in \cite[Section 4, Cases 1--3]{KV}.
Write $N=N_v$ and define, for $T>0$,
\[
 \operatorname{osc}(T)
 =
 \inf_{t_0\in\R}
 \frac{
 \sup_{|t-t_0|\le T N(t_0)^{-2}}N(t)
 }{
 \inf_{|t-t_0|\le T N(t_0)^{-2}}N(t)
 }.
\]

\emph{Case 1: $\operatorname{osc}(T)$
remains bounded as $T\to\infty$.}
We can choose $t_n\in\R$ and $C\ge1$ such that
\[
 C^{-1}\le
 \frac{N(t_n+sN(t_n)^{-2})}{N(t_n)}
 \le C,\quad |s|\le n.
\]
Passing to a normalized time-translate limit
and using the common compactness modulus,
we obtain a nonzero global almost periodic
solution $u$ with \(\inf_{s\in\R}N_u(s)>0\).

\emph{Case 2: $\operatorname{osc}(T)\to\infty$.}
Set
\[
 a_0(t)=
 \frac{N(t)}{\sup_{r\le t}N(r)}
 +\frac{N(t)}{\sup_{r\ge t}N(r)},
\]
where a quotient with infinite denominator
is understood to be zero.

If $\inf_t a_0(t)=0$, the selection in
\cite[Section 4, Case 2]{KV}
gives times $t_n^-<r_n<t_n^+$ such that
\[
 N(r_n)\le
 2\inf_{t\in[t_n^-,t_n^+]}N(t),
 \qquad
 \frac{N(t_n^\pm)}{N(r_n)}\to\infty.
\]
Let $v_n$ be the normalized time translates
at $r_n$, and put \(s_n^\pm=(t_n^\pm-r_n)N(r_n)^2\).
After passing to a subsequence,
$v_n$ converges to a nonzero global
almost periodic solution $u$.

The divergence of the endpoint frequency
ratios and the compactness of the normalized
orbit imply
\[
 s_n^-\to-\infty,\qquad s_n^+\to\infty.
\]
Moreover,
\[
 N_{v_n}(s)
 =\frac{N(r_n+sN(r_n)^{-2})}{N(r_n)}
 \ge\frac12,
 \quad s\in[s_n^-,s_n^+].
\]
By strong convergence on compact time
intervals and the common compactness
modulus, $u$ admits a frequency function
satisfying \(\inf_{s\in\R}N_u(s)>0\).

It remains to exclude $\inf_t a_0(t)>0$.
By the construction in
\cite[Section 4, Case 3]{KV}, after conjugate
time reversal if necessary, there exist
$T>0$ and increasing times $t_n$ such that
\[
 0<t_{n+1}-t_n\le8T N(t_n)^{-2},
 \qquad
 N(t_{n+1})\le\frac12N(t_n).
\]
Normalize at $t_n$ and extract a limit $u$.
By assumption, $u$ is nonzero and global.

The original time $t_0$ has normalized
coordinate \(s_n=(t_0-t_n)N(t_n)^2\).
Since $N(t_n)\le2^{-n}N(t_0)$,
\[
\begin{aligned}
 |s_n|
 \le8T\sum_{k=0}^{n-1}
 \frac{N(t_n)^2}{N(t_k)^2}\le8T\sum_{j=1}^n4^{-j}<8T.
\end{aligned}
\]
However,
\[
 N_{v_n}(s_n)
 =\frac{N(t_0)}{N(t_n)}
 \ge2^n\to\infty.
\]
This contradicts the strong convergence
to a nonzero solution on $[-8T,0]$
and the common compactness modulus.
Thus the last case is impossible, and
one of the first two cases gives a nonzero
global limit $u$ with \(\inf_{t\in\R}N_u(t)>0\).
A fixed critical scaling completes
the proof.
\end{proof}

\begin{corollary}Rigidity under a kinetic gap
\label{cor:tgap-threshold-ap}
There is no nonzero maximal-lifespan strong solution
$v:I\times\R^3\to\C$ to \eqref{eq:nls}
that is almost periodic in $\Hone$ 
and satisfies
\[
E(v)=E(W),\quad
\sup_{t\in I}\|\nabla v(t)\|_{L^2}^2
\le G-\kappa
\]
for some $\kappa>0$.
\end{corollary}

\begin{proof}
Suppose such a solution $v$ exists. By Lemma~\ref{lem:kin-finite-ap}, $v$ is global. By the compactness and stability argument in \cite[Definition~4.1 and the discussion following it]{KV}, every sequence of normalized time translates of $v$ or $\overline{v(-t)}$ has a subsequence converging in $C(J;\Hone)$ on every compact subinterval of the lifespan of the limiting solution. Every such limit is almost periodic, and the two families have a common compactness modulus after adjoining the complex conjugates of the compact closure of the normalized  orbit of $v$. Strong convergence on compact time intervals and invariance under scaling, translation, and conjugate time reversal show that every such limit $u$ satisfies
\[
E(u)=E(W),\qquad
\sup_{t\in I(u)}\|\nabla u(t)\|_{L^2}^2\le G-\kappa.
\]
Thus every such limit is nonzero. By Lemma~\ref{lem:kin-finite-ap}, all these limits are global. Lemma~\ref{lem:closed-class-floor} then gives such a limit $u$ with $\inf_{t\in\R}N_u(t)>0$. After a fixed critical scaling, we may assume $N_u(t)\ge1$ for all $t\in\R$. Proposition~\ref{prop:tgap-forward-rigidity} now implies $u\equiv0$, contradicting $E(u)=E(W)>0$.
\end{proof}

\section{Approach to the ground state orbit}
\label{sec:threshold-compactness}
We first prove global existence and almost periodicity in each nonscattering time direction for threshold solutions on the lower kinetic branch. We then use the rigidity result of Section 3 to show that a forward-nonscattering solution approaches \(\mathcal W\) along a sequence of times.
\subsection{Threshold compactness}

\begin{lemma}
\label{lem:threshold-palais-smale}
Let \(u_n:I_n\times\R^3\to\mathbb C\) be maximal-lifespan solutions,
where \(I_n=(T_n^-,T_n^+)\ni0\), satisfying
\[
 E(u_n)=E(W),\qquad \|\nabla u_n(0)\|_{L^2} \le \|\nabla W\|_{L^2}.
\]
Assume that
\begin{equation}\label{eq:lt-two-sides}
 S_{(T_n^-,0)}(u_n)\to\infty,
 \qquad S_{(0,T_n^+)}(u_n)\to\infty.
\end{equation}
After passing to a subsequence, there are
\(\ell_n>0\), \(x_n\in\R^3\), and \(f\in\Hone\) such that
\[
 g_{\ell_n,x_n}^{-1}u_n(0)\to f
 \quad\text{strongly in }\Hone,
 \qquad E(f)=E(W),\qquad 
 \|\nabla f\|_{L^2} \le \|\nabla W\|_{L^2}.
\]
The maximal-lifespan solution through \(f\) scatters in neither
time direction.
\end{lemma}

\begin{proof}
We follow the profile decomposition argument of
\cite[Lemma~2.5]{DM} and
\cite[Proposition~3.1]{KV}.

Applying \eqref{eq:prelim-profile} to $u_n(0)$,
kinetic and energy decoupling give
\begin{equation}\label{eq:lt-energy-decoupling}
E(W)=\sum_{j=1}^J E(U^j)+E(w_n^J)+o_n(1).
\end{equation}
In particular,
$\|\nabla\phi^j\|_2^2\le G$.
If $\sigma_n^j=0$, the sharp Sobolev
inequality gives
$E(U^j)\ge\frac13\|\nabla\phi^j\|_2^2$.
If $\sigma_n^j\to\pm\infty$, then
$E(U^j)=\frac12\|\nabla\phi^j\|_2^2$
by \eqref{eq:free-L6-decay}.
Thus, in either case,
\begin{equation}\label{eq:lt-profile-positivity}
E(U^j)\ge
\frac13\|\nabla\phi^j\|_2^2\ge0.
\end{equation}
Similarly, for each fixed $J$ and all
large $n$,
\begin{equation}\label{eq:lt-remainder-positivity}
E(w_n^J)\ge
\frac14\|\nabla w_n^J\|_2^2.
\end{equation}

Suppose every profile has energy strictly
below $E(W)$. By
\eqref{eq:lt-energy-decoupling} and
\eqref{eq:lt-profile-positivity}, these
energies have a uniform gap below $E(W)$.
We verify the kinetic condition needed
for the below-ground state scattering theorem.

If $\sigma_n^j=0$, then
$\|\nabla\phi^j\|_2^2<G$; otherwise
the sharp Sobolev inequality would give
$E(U^j)\ge E(W)$.
If $\sigma_n^j\to\pm\infty$, then
\[
\|\nabla\phi^j\|_2^2
=2E(U^j)<\frac{2G}{3}<G.
\]
Since $U^j$ scatters in the corresponding
direction, its kinetic norm is strictly
below $G$ at sufficiently large positive
or negative times.

Thus every nonlinear profile satisfies
the hypotheses of the below ground state
scattering theorem of \cite{SZ} and is
global and scattering. By kinetic
decoupling, all sufficiently high-index
profiles are controlled by the small-data
theory. Nonlinear profile approximation
and stability then give uniform scattering
bounds for $u_n$, contradicting
\eqref{eq:lt-two-sides}.

Hence some profile has energy $E(W)$.
By \eqref{eq:lt-energy-decoupling}--%
\eqref{eq:lt-remainder-positivity},
all other profiles vanish and
$\|\nabla w_n^1\|_2\to0$.
Consequently, after relabeling,
\begin{equation}\label{eq:lt-single-linear-profile}
u_n(0)=g_{\ell_n,x_n}
e^{i\sigma_n\Delta}\phi+o_n(1)
\quad\text{in }\Hone.
\end{equation}

If $\sigma_n\to+\infty$ or
$\sigma_n\to-\infty$, the free Strichartz
norm tends to zero in the corresponding
time direction. The small-data theory
then contradicts \eqref{eq:lt-two-sides}.
Thus $\sigma_n=0$, and
\[
g_{\ell_n,x_n}^{-1}u_n(0)
\to f:=\phi
\quad\text{strongly in }\Hone.
\]
By continuity and scaling invariance,
\[
E(f)=E(W),\quad
\|\nabla f\|_2\le\|\nabla W\|_2.
\]

Finally, if the solution with initial
data $f$ scattered in either time
direction, stability would give a
uniform scattering bound for $u_n$
in that direction, again contradicting
\eqref{eq:lt-two-sides}.
\end{proof}

\begin{proposition}
\label{prop:lt-global-ap}
Let $u:I\times\R^3\to\C$ be a
maximal-lifespan strong solution to \eqref{eq:nls}
with $0\in I$, satisfying
\[
 E(u)=E(W),\quad
 \|\nabla u_0\|_{L^2}<\|\nabla W\|_{L^2}.
\]
Then $I=\R$.
Moreover, if $u$ does not scatter forward in time,
then $u$ is almost periodic on $[0,\infty)$
modulo scaling and translation.
The analogous conclusion holds backward in time.
\end{proposition}

\begin{proof}
By the variational properties in
Section~\ref{sec:preliminaries},
\[
 \|\nabla u(t)\|_{L^2}
 <\|\nabla W\|_{L^2},
 \quad t\in I.
\]

If $S_{[0,T_+)}(u)<\infty$, the scattering
criterion gives $T_+=\infty$.
Otherwise, let $t_n\to T_+$ with $t_n>0$.
Then
\[
 S_{(T_-,t_n)}(u)\to\infty,
 \quad
 S_{(t_n,T_+)}(u)=\infty.
\]
Applying Lemma~\ref{lem:threshold-palais-smale}
to the time translates $u(t_n+\cdot)$,
we obtain, after passing to a subsequence,
\[
 g_{\ell_n,x_n}^{-1}u(t_n)\to f
 \quad\text{strongly in }\Hone.
\]
Since the sequence $t_n$ was arbitrary
and $u$ is continuous in $\Hone$,
the forward orbit is precompact modulo
scaling and translation. 
The forward-endpoint argument of
Lemma~\ref{lem:kin-finite-ap}
then gives $T_+=\infty$.
Applying the same argument to
$\overline{u(-t)}$ yields $T_-=-\infty$
and the corresponding backward
almost periodicity whenever $u$
does not scatter backward.
\end{proof}

\subsection{Sequential approach}

Corollary~\ref{cor:tgap-threshold-ap} excludes a nonzero global
almost periodic solution with $E=E(W)$ and a uniform kinetic gap.
The next lemma puts a hypothetical forward compact tail into this
whole-lifespan setting.

\begin{lemma}
\label{lem:lt-gap-tail-limit}
Let \(u\) be a forward-nonscattering solution as in
Proposition~\ref{prop:lt-global-ap}. Suppose that, for some \(T_0\ge0\)
and \(\kappa>0\),
\[
 \sup_{t\ge T_0} \|\nabla u(t)\|_{L^2}^2\le G-\kappa.
\]
Then there exists a nonzero global almost periodic solution \(v\)
with \(E(v)=E(W)\) and
\[
 \sup_{t\in\R} \|\nabla v(t)\|_{L^2}^2\le G-\kappa.
\]
\end{lemma}

\begin{proof}
Let $t_n\to\infty$. Since $u$ does not
scatter forward,
\[
 S_{(-\infty,t_n)}(u)\to\infty,
 \quad S_{(t_n,\infty)}(u)=\infty.
\]
By Lemma~\ref{lem:threshold-palais-smale},
after passing to a subsequence, there are
$\ell_n>0$, $x_n\in\R^3$, and $v_0\in\Hone$
such that
\[
 g_{\ell_n,x_n}^{-1}u(t_n)\to v_0
 \quad\text{strongly in }\Hone.
\]
Moreover,
\[
 E(v_0)=E(W),\quad
 \|\nabla v_0\|_{L^2}^2\le G-\kappa,
\]
and the solution $v$ with $v(0)=v_0$
scatters in neither time direction.
Proposition~\ref{prop:lt-global-ap}
therefore implies that $v$ is global
and almost periodic on $\R$.
Define
\[
 v_n(s)=g_{\ell_n,x_n}^{-1}
 u(t_n+\ell_n^2s).
\]
By stability, for every compact
$J\subset\R$,
\[
 v_n\to v\quad\text{in }C(J;\Hone),
 \qquad \sup_n S_J(v_n)<\infty.
\]
We claim that \((t_n-T_0)/{\ell_n^2} \to\infty\).
Otherwise, along a subsequence,
$(t_n-T_0)/\ell_n^2\le C$ for some $C>0$.
Hence
\[
\begin{aligned}
 S_{[T_0,t_n]}(u)
 =S_{[-(t_n-T_0)/\ell_n^2,0]}(v_n)\le S_{[-C,0]}(v_n)\le C',
\end{aligned}
\]
contrary to the forward nonscattering
assumption.
Thus, for every fixed $s\in\R$, \(t_n+\ell_n^2s\ge T_0\)
for all sufficiently large $n$.
Strong convergence and scaling invariance
then give
\[
\begin{aligned}
 \|\nabla v(s)\|_{L^2}^2
 =\lim_{n\to\infty}
 \|\nabla u(t_n+\ell_n^2s)\|_{L^2}^2\le G-\kappa.
\end{aligned}
\]
Since $s\in\R$ is arbitrary, the desired
kinetic bound follows. Finally,
$E(v)=E(W)>0$ implies $v\not\equiv0$.
\end{proof}

\begin{corollary}
\label{cor:lt-sequential-ground}
If \(u\) satisfies
\[
 E(u)=E(W),\qquad \|\nabla u_0\|_{L^2}<\|\nabla W\|_{L^2}
\]
and does not scatter forward, then for every \(T\ge0\),
\begin{equation}\label{eq:lt-tail-sup}
 \sup_{t\ge T} \|\nabla u(t)\|_{L^2}^2=G.
\end{equation}
Consequently there are times \(t_n\to\infty\) for which
\begin{equation}\label{eq:lt-sequential-ground}
 \|\nabla u(t_n)\|_{L^2}\to \|\nabla W\|_{L^2},\qquad d_{\Wcal}(u(t_n))\to0.
\end{equation}
\end{corollary}

\begin{proof}
By Proposition~\ref{prop:lt-global-ap},
$u$ is global and forward almost periodic.
Suppose \eqref{eq:lt-tail-sup} fails.
Since $\|\nabla u(t)\|_{L^2}^2<G$
for all $t\in\R$, there exist
$T_0\ge0$ and $\kappa>0$ such that \(\sup_{t\ge T_0}\|\nabla u(t)\|_{L^2}^2
 \le G-\kappa\).
By Lemma~\ref{lem:lt-gap-tail-limit},
there exists a nonzero global almost
periodic solution satisfying the conditions
excluded by Corollary~\ref{cor:tgap-threshold-ap}.
This is a contradiction, proving
\eqref{eq:lt-tail-sup}.

Choose $t_n\ge n$ such that
\[
 \|\nabla u(t_n)\|_{L^2}^2>G-\frac1n.
\]
By energy conservation,
\[
 \|u(t_n)\|_{L^6}^6
 =3\|\nabla u(t_n)\|_{L^2}^2-2G
 \to G.
\]
By forward almost periodicity, after passing
to a subsequence and applying suitable
scaling and translation, $u(t_n)$ converges
strongly in $\Hone$ to some $f$ satisfying \(\|\nabla f\|_{L^2}^2
 =\|f\|_{L^6}^6=G\).
Equality in the sharp Sobolev inequality
implies $f\in\Wcal$.
Hence, after relabeling the subsequence,
\eqref{eq:lt-sequential-ground} follows.
\end{proof}

\section{Dynamics near the ground state orbit}
\label{sec:formal-local}

To obtain classification, we must rule out later exits and identify solutions that remain near \(\mathcal W\). This section develops the local estimates needed for these two tasks. Modulation controls the symmetry parameters, while comparison with the special solutions identifies exit limits and exponentially convergent solutions.

\subsection{Spectral decomposition and coercivity}

We begin with the linearization around the ground state $W$. The aim is to separate the two real eigenmodes and the directions generated by the symmetries of $W$, and then establish coercivity on the remaining space.
Write $u=W+h$, where $h=h_1+ih_2$, and set
\[
L_+=-\Delta-5W^4,\quad L_-=-\Delta-W^4,\quad
\mathcal A=\begin{pmatrix}0&L_-\\-L_+&0\end{pmatrix}.
\]
Then
\begin{equation}\label{eq:exit-fixed-equation}
 \partial_t h=\mathcal Ah+\Ncal(h),\quad \Ncal(h)=i\left[
 |W+h|^4(W+h)-W^5
 -W^4(3h+2\overline h)
 \right].
\end{equation}
We identify complex-valued functions with pairs of real functions. The quadratic part of the energy at $W$ is
\[
Q(h)=\frac12\int_{\R^3}
\bigl(|\nabla h|^2-W^4(5h_1^2+h_2^2)\bigr)\dd x,
\]
with associated real bilinear form
\[
Q(f,g)=\frac12\bigl[Q(f+g)-Q(f)-Q(g)\bigr].
\]

By \cite[Lemma 5.1 and (5.7)]{DM}, there are radial functions $Y_\pm\in\mathcal S(\R^3)$ and $e_0>0$ such that
\[
\mathcal AY_\pm=\pm e_0Y_\pm,\qquad Y_-=\overline{Y_+}.
\]
Writing $Y_+=Y_1+iY_2$, the eigenvalue equations give
\[
Q(Y_\pm)=0,\qquad
Q(Y_+,Y_-)=-\langle L_-Y_2,Y_2\rangle<0,
\]
where the positivity of $L_-$ and its kernel are established below. Set
\[
q_*=-2Q(Y_+,Y_-)>0.
\]

The symmetries of $W$ give the zero modes
\[
Z=\operatorname{span}_{\R}
\{iW,\Lambda W,\partial_1W,\partial_2W,\partial_3W\},
\quad \Lambda W=\frac12W+x\cdot\nabla W.
\]
In particular, $Q(z,h)=0$ for every $z\in Z$ and $h\in\Hone$.
Define
\[
a_\pm(h)=\frac{Q(h,Y_\mp)}{Q(Y_+,Y_-)},
\quad P_\pm h=a_\pm(h)Y_\pm,\quad
P_c=I-P_+-P_-.
\]
Let $\Pi_Z$ be the $\Hone$-orthogonal
projection onto $Z$. Since $Q(Z,\Hone)=0$,
we have $Z\subset P_c\Hone$. Thus every
$h\in\Hone$ has the decomposition
\[
h=a_+(h)Y_++a_-(h)Y_-
+\Pi_ZP_ch+(I-\Pi_Z)P_ch,
\]
where
\[
 (I-\Pi_Z)P_ch\in\mathcal H_\perp,
 \qquad
 \mathcal H_\perp
 :=
 P_c\Hone\cap Z^{\perp_{\Hone}}.
\]
Moreover,
\[
Q(h)=-q_*a_+(h)a_-(h)
+Q((I-\Pi_Z)P_ch).
\]
The projections are bounded on $\Hone$. It remains to prove the nonradial spectral properties and the coercivity of $Q$ on $\mathcal H_\perp$.

\begin{lemma}
\label{lem:local-sphere-form}
The quadratic form $Q$ has negative index one, and
its bilinear kernel is
\[
 Z=\operatorname{span}_{\mathbb R}
 \{iW,\Lambda W,\partial_1W,\partial_2W,\partial_3W\}.
\]
\end{lemma}

\begin{proof}
Set $y=x/\sqrt{3}$ and define
\[
 \Omega(y)=\frac{2}{1+|y|^2},
 \quad
 \omega(y)=\frac{(2y,|y|^2-1)}{1+|y|^2}.
\]
The map $\omega$ identifies $\mathbb R^3$ with
$\mathbb S^3$ minus the north pole, and the pullback
of the round metric is $\Omega^2dy^2$.
Moreover,
\[
 W(\sqrt{3}y)=\frac{\Omega(y)^{1/2}}{\sqrt{2}}.
\]

For $f=f_1+if_2\in H^1(\mathbb S^3)$, let \( h(x)=W(x)f\left(\omega\left(\frac{x}{\sqrt{3}}\right)\right)\).
The conformal Laplacian identity gives
\[
 -\Delta_y\bigl(\Omega^{1/2}f\circ\omega\bigr)
 =
 \Omega^{5/2}
 \left(-\Delta_{\mathbb S^3}f+\frac34f\right)
 \circ\omega.
\]
Using $x=\sqrt{3}y$ and integrating by parts, we obtain
\begin{equation}
 \|h\|_{\dot H^1}^2
 =
 \frac{\sqrt{3}}{2}
 \int_{\mathbb S^3}
 \left(
 |\nabla_{\mathbb S^3}f|^2+\frac34|f|^2
 \right)\,d\omega.
 \label{eq:local-sphere-norm}
\end{equation}
Since $d\omega=\Omega^3dy$, we also have
\[
 \int_{\mathbb R^3}W^4|h|^2\,dx
 =
 \frac{3\sqrt{3}}8
 \int_{\mathbb S^3}|f|^2\,d\omega.
\]
Consequently,
\begin{equation}
 Q(h)
 =
 \frac{\sqrt{3}}4
 \int_{\mathbb S^3}
 \left(
 |\nabla_{\mathbb S^3}f_1|^2-3f_1^2
 +|\nabla_{\mathbb S^3}f_2|^2
 \right)\,d\omega.
 \label{eq:local-sphere-Q}
\end{equation}
These identities first hold for smooth functions
supported away from the north pole. Such functions
are dense in $H^1(\mathbb S^3)$, since a cutoff
around the north pole on a ball of radius
$\varepsilon$ has squared gradient integral
$O(\varepsilon)$.
By \eqref{eq:local-sphere-norm} and density,
the transformation is an isomorphism from
$H^1(\mathbb S^3)$ onto $\dot H^1(\mathbb R^3)$.

The eigenvalues of $-\Delta_{\mathbb S^3}$ are \(\ell(\ell+2)\), \(\ell=0,1,\ldots\).
By \eqref{eq:local-sphere-Q}, the real part has
eigenvalues $\ell(\ell+2)-3$, which are negative
only for $\ell=0$ and vanish only for $\ell=1$.
The imaginary part is nonnegative and has
kernel consisting of the constants.
Thus $Q$ has exactly one negative direction.
Its real kernel is spanned by the four first
spherical harmonics, while its imaginary
kernel is spanned by the constants.

Finally,
\[
 \frac{\Lambda W(x)}{W(x)}
 =-\frac12\omega_4\left(\frac{x}{\sqrt{3}}\right),
 \quad
 \frac{\partial_jW(x)}{W(x)}
 =-\frac{1}{2\sqrt{3}}
 \omega_j\left(\frac{x}{\sqrt{3}}\right).
\]
These identities identify the real kernel with
$\operatorname{span}_{\mathbb R}
\{\Lambda W,\partial_1W,\partial_2W,\partial_3W\}$,
and the imaginary kernel with
$\operatorname{span}_{\mathbb R}\{iW\}$.
Hence the bilinear kernel of $Q$ is $Z$.
\end{proof}

\begin{lemma}
\label{lem:exit-linear}
The projections \(P_{\pm}\),\(P_c\) commute with \(e^{t\mathcal A}\). On \(X_c=P_c\Hone\),
\begin{equation}\label{eq:exit-center-coercivity}
 Q(h)\ge c\,\dist_{\Hone}(h,Z)^2,
 \qquad
 \|e^{t\mathcal A}P_ch\|_{\Hone}
       \le C(1+|t|)\|P_ch\|_{\Hone}.
\end{equation}
For an interval \(J\) with \(|J|\le 1\), use $X(J)$ from
\eqref{eq:prelim-X} and set
\[
 \|F\|_{N(J)}
   =\|F\|_{L_t^2L_x^{\frac{6}{5}}(J)}
        +\|\nabla F\|_{L_t^2L_x^{\frac{6}{5}}(J)} .
\]
If \(\partial_t h=\mathcal A h+F\) on \(J\), then
\begin{equation}\label{eq:exit-local-linear}
 \|h\|_{X(J)}
       \le C\bigl(\|h(t_0)\|_{\Hone}+\|F\|_{N(J)}\bigr),
 \qquad t_0\in\partial J .
\end{equation}
For \(\|h\|_{X(J)}+\|k\|_{X(J)}\) sufficiently small,
\begin{equation}\label{eq:exit-nonlinear}
 \|\Ncal(h)-\Ncal(k)\|_{N(J)}
 \le C\bigl(\|h\|_{X(J)}+\|k\|_{X(J)}\bigr)
                       \|h-k\|_{X(J)} .
\end{equation}
\end{lemma}

\begin{proof}
We regard $\Hone$ as a real Hilbert space
and write $\mathcal A=\mathcal A_0+B$, where
\[
\mathcal A_0h=i\Delta h,\quad
Bh=iW^4(3h+2\overline h).
\]
The free Schr\"odinger group
$e^{t\mathcal A_0}=e^{it\Delta}$ is a
strongly continuous unitary group on $\Hone$.
Moreover, since $W^4\in L^\infty$ and
$\nabla(W^4)\in L^3$, Sobolev embedding gives
\[
\begin{aligned}
\|Bh\|_{\Hone}
\lesssim
\|W^4\|_\infty\|\nabla h\|_2
+\|\nabla(W^4)\|_3\|h\|_6\lesssim\|h\|_{\Hone}.
\end{aligned}
\]
Thus $B$ is bounded and real-linear on
$\Hone$. By the bounded perturbation
theorem for strongly continuous groups,
$\mathcal A$ generates a $C_0$-group
$e^{t\mathcal A}$ on $\Hone$.

For data in the domain of $\mathcal A$,
the identity
$Q(\mathcal Af,g)=-Q(f,\mathcal Ag)$
and the eigenvalue equations give
\[
Q(e^{t\mathcal A}h)=Q(h),\quad
a_\pm(e^{t\mathcal A}h)
=e^{\pm e_0t}a_\pm(h).
\]
Thus $P_\pm$ and $P_c$ commute with
$e^{t\mathcal A}$. By density and
strong continuity, these identities
extend to all data in $\Hone$.

We next prove coercivity.
By Lemma~\ref{lem:local-sphere-form},
$Q$ has negative index one and bilinear
kernel $Z$. Its restriction to
$\operatorname{span}_{\R}\{Y_+,Y_-\}$
is nondegenerate with one negative
direction. Since $X_c$ is the
$Q$-orthogonal complement of this
space, $Q$ is nonnegative on $X_c$.

Suppose the coercivity estimate fails.
There exist
$h_n\in X_c\cap Z^{\perp_{\Hone}}$
such that
\[
\|h_n\|_{\Hone}=1,\quad Q(h_n)\to0.
\]
Passing to a subsequence, assume
$h_n\rightharpoonup h$ weakly in $\Hone$.
Set
\[
V(v)=\int_{\R^3}W^4
\left(5(\Repart v)^2+(\Impart v)^2\right)
\dd x.
\]
Since $Q(h_n)\to0$, we have
$V(h_n)\to1$. As $W^4\in L^{3/2}$,
local Rellich compactness and Sobolev
embedding imply $V(h_n)\to V(h)$.
Hence $V(h)=1$.

The space $X_c\cap Z^{\perp_{\Hone}}$
is weakly closed, and weak lower
semicontinuity gives
$\|h\|_{\Hone}\le1$. Consequently,
\[
0\le Q(h)
=\frac12\|h\|_{\Hone}^2-\frac12
\le0.
\]
Thus $Q(h)=0$. Since $Q$ is nonnegative
on $X_c$, polarization gives
$Q(h,X_c)=0$. Together with the
$Q$-orthogonality of $X_c$ and
$\operatorname{span}_{\R}\{Y_+,Y_-\}$,
this implies $h\in Z$.
But $h\perp_{\Hone}Z$, so $h=0$,
contradicting $V(h)=1$.
Therefore,
\[
Q(h)\ge
c\|h-\Pi_Zh\|_{\Hone}^2,
\quad h\in X_c.
\]

For $h\in X_c$, write
\[
e^{t\mathcal A}h=g(t)+z(t),
\quad g(t)\perp_{\Hone}Z,\quad z(t)\in Z.
\]
Since $Z\subset X_c$, both $g(t)$
and $z(t)$ belong to $X_c$.
Using $Q(Z,\Hone)=0$, conservation
of $Q$, and the coercivity just proved,
we obtain
\[
\begin{aligned}
\|g(t)\|_{\Hone}^2
\le C Q(g(t))
=CQ(h)\le C\|h\|_{\Hone}^2.
\end{aligned}
\]

Let $\{Z_j\}_{j=0}^4$ be a basis of $Z$,
and let $\{\psi_j\}_{j=0}^4\subset Z$
be its dual basis for the real
$\Hone$ inner product. Set
$\ell_j(v)=\langle v,\psi_j\rangle_{\Hone}$.
For smooth $v$, integration by parts
and the complex representation of
$\mathcal A$ give
\[
\begin{aligned}
|\ell_j(\mathcal Av)|
\lesssim\|v\|_{L^6}
\left(
\|\Delta^2\psi_j\|_{L^{6/5}}
+\|W^4\Delta\psi_j\|_{L^{6/5}}
\right)\lesssim\|v\|_{\Hone}.
\end{aligned}
\]
The required integrability follows
from the decay of $W$ and its derivatives.
Thus these coefficient pairings extend
continuously to $\Hone$.

Writing $z(t)=\sum_{j=0}^4z_j(t)Z_j$,
we have $z_j(t)=\ell_j(e^{t\mathcal A}h)$.
Since $\mathcal AZ=0$, for data in
the domain of $\mathcal A$,
\[
z_j'(t)=\ell_j(\mathcal Ag(t)),\quad
|z_j'(t)|\le C\|g(t)\|_{\Hone}
\le C\|h\|_{\Hone}.
\]
Integrating in time gives
\[
\|z(t)\|_{\Hone}
\le C(1+|t|)\|h\|_{\Hone}.
\]
Together with the bound for $g(t)$,
this proves the second estimate in
\eqref{eq:exit-center-coercivity}.
The general case follows by density.

To prove \eqref{eq:exit-local-linear},
write the linearized equation as
\[
\partial_t h=i\Delta h
+iW^4(3h+2\overline h)+F.
\]
Since $W^4\in L^{3/2}\cap L^3$
and $\nabla(W^4)\in L^{3/2}$,
H\"older's inequality and Sobolev
embedding give
\[
\|W^4h\|_{N(J)}
\le C|J|^{1/2}
\|h\|_{L_t^\infty(J;\Hone)}.
\]
Free Strichartz estimates yield
\[
\begin{aligned}
\|h\|_{X(J)}
\le C\bigl(
\|h(t_0)\|_{\Hone}+\|F\|_{N(J)}
\bigr)+C|J|^{1/2}\|h\|_{X(J)}.
\end{aligned}
\]
Choosing $|J|$ sufficiently small
allows the last term to be absorbed.
A finite subdivision gives
\eqref{eq:exit-local-linear}
for $|J|\le1$.

Finally, $\Ncal(h)$ is a sum of terms
$W^{5-m}h^a\overline h^{\,b}$, where
$a+b=m$ and $2\le m\le5$.
Since $W,\nabla W\in L^6$, each
differentiated product is estimated
using H\"older's inequality with five
spatial $L^6$ factors. For example,
\[
\begin{aligned}
\|W^3(h+k)\nabla(h-k)\|_{L_t^2L_x^{6/5}}\lesssim
\|W\|_{L^6}^3
\left(
\|h\|_{L_t^\infty L_x^6}
+\|k\|_{L_t^\infty L_x^6}
\right)
\|\nabla(h-k)\|_{L_t^2L_x^6}.
\end{aligned}
\]
Terms with a derivative on $W$
gain a factor $|J|^{1/2}$.
The remaining terms are estimated
similarly. Thus, with
$M=\|h\|_{X(J)}+\|k\|_{X(J)}$,
\[
\|\Ncal(h)-\Ncal(k)\|_{N(J)}
\le C\sum_{m=2}^5M^{m-1}
\|h-k\|_{X(J)}.
\]
For $M$ sufficiently small,
this proves \eqref{eq:exit-nonlinear}.
\end{proof}

\subsection{Modulation estimates}
\label{subsec:audited-threshold-exit}

\begin{lemma}[Modulation near the ground state orbit]
\label{lem:general-first-exit-coercivity}
There are $\delta_0>0$ and $C<\infty$ with the following properties.
For every $f$ with $d_{\mathcal W}(f)<\delta_0$, there are phase,
frequency scale and center parameters for which
\[
 f=\mathcal T_{\theta,\lambda,a}(W+v),\qquad
 v\perp_{\dot H^1}Z,
 \qquad
 \|v\|_{\dot H^1}\sim d_{\mathcal W}(f).
\]
The parameters are smooth local functionals of $f$ in the energy
topology, with the phase taken modulo $2\pi$.
Assume in addition that $E(f)=E(W)$, and set
\[
 \xi_1=\frac{a_+(v)+a_-(v)}{2},\quad
 \xi_2=\frac{a_+(v)-a_-(v)}{2},\quad \gamma=P_cv.
\]
Then
\begin{align}
 q_*\xi_2^2+Q(\gamma)
 &=q_*\xi_1^2+R_E(v),
 &|R_E(v)|&\le C\|v\|_{\dot H^1}^3,
 \label{eq:formal-mod-energy}\\
 \|v\|_{\dot H^1}&\le C|\xi_1|,
 &Q(\gamma)&\ge0.
 \label{eq:general-scalar-coercivity}
\end{align}
Along a strong solution of energy $E(W)$ staying in this tube, the
parameters can be chosen continuously, with a continuous real lift of the phase.
They and the scalar coefficients are continuously differentiable
in the rescaled time $\frac{\mathrm{d}\tau}{\mathrm{d}t}=\lambda^2$, and satisfy
\begin{equation}\label{eq:general-scalar-normal-form}
    (\xi_1)_\tau=e_0\xi_2+O(\|v\|_{\dot H^1}^2),
    \qquad
 (\xi_2)_\tau=e_0\xi_1+O(\|v\|_{\dot H^1}^2),
\end{equation}
\begin{equation}\label{eq:formal-mod-speed}
    |\theta_\tau|+|\partial_\tau\log\lambda|+|\lambda a_\tau|
 \le C\|v\|_{\dot H^1}.
\end{equation}
In particular, for a nonstationary solution with $E(u)=E(W)$,
$\xi_1$ has a constant nonzero sign on each connected tube interval.
\end{lemma}

\begin{proof}
Let
\[
Z_0=iW,\quad Z_1=\Lambda W,\quad
Z_{j+1}=\partial_jW,\quad 1\le j\le3,
\]
and write $q=(\theta,\log\lambda,a)$.
Define \(\mathcal F_k(f,q)
=
\langle\mathcal T_q^{-1}f-W,Z_k
\rangle_{\Hone}\).
Since $\mathcal T_q$ is an isometry of
$\Hone$,
\[
\mathcal F_k(f,q)
=
\langle f,\mathcal T_qZ_k\rangle_{\Hone}
-\langle W,Z_k\rangle_{\Hone}.
\]
The smoothness and decay of $Z_k$ imply
that $\mathcal F$ is smooth in $(f,q)$.
At $(W,0)$, its derivative in $q$ is,
up to column signs, the Gram matrix of
$Z_0,\ldots,Z_4$, which is invertible.
The implicit function theorem gives
locally unique parameters such that
\[
f=\mathcal T_q(W+v),\quad
v\perp_{\Hone}Z,\quad
\|v\|_{\Hone}\lesssim\|f-W\|_{\Hone}.
\]

Applying this construction after undoing
a symmetry of a ground state within
$2d_{\mathcal W}(f)$ gives
\[
\|v\|_{\Hone}\sim d_{\mathcal W}(f)
\]
uniformly near the full orbit.
The parameters patch continuously by
local uniqueness and the properness of
the symmetry action on $W$.
Indeed, strong convergence
$\mathcal T_{q_n}W\to W$ forces
$q_n\to0$, modulo $2\pi$ in phase.

Suppose $E(f)=E(W)$. The spectral
decomposition gives
\[
v=(\xi_1+\xi_2)Y_+
+(\xi_1-\xi_2)Y_-+\gamma,
\quad \gamma=P_cv.
\]
Since $v\perp_{\Hone}Z$, we have
$\|\Pi_Z\gamma\|_{\Hone}
\lesssim|\xi_1|+|\xi_2|$.
Lemma~\ref{lem:exit-linear} therefore gives
\[
\begin{aligned}
Q(\gamma)&\ge
c\|\gamma-\Pi_Z\gamma\|_{\Hone}^2,\\
\|v\|_{\Hone}^2
&\lesssim\xi_1^2+\xi_2^2+Q(\gamma).
\end{aligned}
\]

Since $E'(W)=0$, Taylor expansion and
the spectral orthogonality relations give
\[
E(W+v)-E(W)=Q(v)+R_3(v),
\]
\[
Q(v)=q_*(\xi_2^2-\xi_1^2)+Q(\gamma)
\]
where $|R_3(v)|\lesssim\|v\|_{\Hone}^3$.
The energy constraint proves
\eqref{eq:formal-mod-energy} with
$R_E=-R_3$. Substitution yields
\[
\|v\|_{\Hone}^2
\lesssim\xi_1^2+\|v\|_{\Hone}^3.
\]
Taking $\delta_0$ sufficiently small,
we absorb the cubic term and obtain
$\|v\|_{\Hone}\lesssim|\xi_1|$.
This proves
\eqref{eq:general-scalar-coercivity}.

We next consider a strong solution
remaining in the tube. Since
$u_t\in C_t\dot H^{-1}$ and
$-\Delta(\mathcal T_qZ_k)$ depends
smoothly on $q$ in $\Hone$,
\[
\partial_t\mathcal F_k(u(t),q)
=
\langle u_t(t),-\Delta(\mathcal T_qZ_k)
\rangle_{\dot H^{-1},\Hone}.
\]
Thus the orthogonality equations are
$C^1$ in $(t,q)$. The implicit function
theorem gives $C^1$ modulation parameters,
also in the rescaled time $\tau$.

In moving coordinates,
\begin{equation}\label{eq:mod-moving-equation}
\begin{aligned}
v_\tau=\mathcal Av+\Ncal(v)
-\theta_\tau i(W+v)-(\log\lambda)_\tau\Lambda(W+v)
+\lambda a_\tau\cdot\nabla(W+v)
\end{aligned}
\end{equation}
in the distributional sense.

Write
$\beta=(-\theta_\tau,
-(\log\lambda)_\tau,\lambda a_\tau)$.
The generators $i,\Lambda,\partial_j$
are skew-adjoint for the real $\Hone$
inner product. Their pairings with
$v\in\Hone$ are therefore defined by
moving the generators onto the fixed
functions $Z_k$, and are bounded by
$C\|v\|_{\Hone}$.
Differentiating the orthogonality
conditions gives
\[
\begin{aligned}
\sum_{j=0}^4
\langle Z_j,Z_k\rangle_{\Hone}\beta_j
=-\langle\mathcal Av+\Ncal(v),
-\Delta Z_k\rangle_{\dot H^{-1},\Hone}+O(|\beta|\|v\|_{\Hone}).
\end{aligned}
\]
Here $-\Delta Z_k\in\Hone$, while
$\mathcal A$ maps $\Hone$ boundedly into
$\dot H^{-1}$. Moreover,
$\|\Ncal(v)\|_{L^{6/5}}
\lesssim\|v\|_{\Hone}^2$
for small $v$. Thus the dual pairing
on the right is $O(\|v\|_{\Hone})$.
Since the Gram matrix is invertible,
we may absorb the error for small
$\delta_0$ and obtain
\[
|\theta_\tau|+|(\log\lambda)_\tau|
+|\lambda a_\tau|
\lesssim\|v\|_{\Hone}.
\]
This proves \eqref{eq:formal-mod-speed}.

Finally, we derive the spectral equations.
The functionals $a_\pm$ are defined by
pairing with $Y_\mp$ through $Q$.
Moving the elliptic operators onto
these fixed smooth functions extends
the pairings to the distributional
terms in \eqref{eq:mod-moving-equation}.
The identity
$Q(\mathcal Af,g)=-Q(f,\mathcal Ag)$
and the eigenvalue equations give
\[
a_\pm(\mathcal Av)=\pm e_0a_\pm(v),
\quad a_\pm(Z_j)=0.
\]
Moving the symmetry generators onto
the fixed test functions also gives
\[
|a_\pm(iv)|+|a_\pm(\Lambda v)|
+\sum_{j=1}^3|a_\pm(\partial_jv)|
\lesssim\|v\|_{\Hone},
\]
\[
|a_\pm(\Ncal(v))|
\lesssim\|v\|_{\Hone}^2.
\]
Applying $a_\pm$ to
\eqref{eq:mod-moving-equation} and using
\eqref{eq:formal-mod-speed}, we obtain
\[
\partial_\tau a_\pm(v)
=\pm e_0a_\pm(v)
+O(\|v\|_{\Hone}^2).
\]
Taking the sum and difference proves
\eqref{eq:general-scalar-normal-form}.
The fixed-test representation also
gives continuity of these derivatives.

If $\xi_1=0$ at some time, then
\eqref{eq:general-scalar-coercivity}
gives $v=0$. Uniqueness of the NLS flow
then forces $u$ to be stationary.
Thus, for a nonstationary solution,
$\xi_1$ has a constant nonzero sign
on every connected tube interval.
\end{proof}

The scalar equations give an integrable bound for the distance to
the ground state orbit.

\begin{lemma}[An integrated modulation estimate]
\label{lem:local-integrated-distance}
Let $u:I\times\R^3\to\C$ be a nonstationary
strong solution to \eqref{eq:nls} 
on an interval $I$ with
$E(u)=E(W)$ and $d_{\mathcal W}(u(t))\le\delta$, where
$\delta>0$ is sufficiently small.
Let $\xi_1,\xi_2$ and $\tau$ be as in
Lemma~\ref{lem:general-first-exit-coercivity}.
Set
\[
\sigma=\operatorname{sgn}\xi_1,\quad
b=\sigma\xi_1,\quad c=\sigma\xi_2.
\]
Then $\sigma$ is constant, $b>0$, and
\begin{equation}\label{eq:local-b-size}
\|v\|_{\Hone}\sim d_{\Wcal}(u)\sim b,
\quad b+|c|\le C\delta.
\end{equation}
Moreover, with primes denoting
differentiation in $\tau$,
\begin{equation}\label{eq:local-bc-normal}
b'=e_0c+O(b^2),\quad
c'=e_0b+O(b^2),
\end{equation}
and
\begin{equation}\label{eq:local-monotone-c}
c'\ge\frac{e_0}{2}b>0,\quad
|c|+|b'|\le Cb.
\end{equation}
Consequently, for every compact interval
$J\subset\tau(I)$,
\begin{equation}\label{eq:local-distance-L1}
\int_J b(\tau)\dd\tau\le C\delta,
\end{equation}
where $C$ is independent of $J$.
\end{lemma}

\begin{proof}
By Lemma~\ref{lem:general-first-exit-coercivity} and the boundedness of \(a_{\pm}\),
$\xi_1$ has a constant nonzero sign, and
\[
\|v\|_{\Hone}\sim d_{\Wcal}(u)\sim|\xi_1|,
\quad
|\xi_2|\le C\|v\|_{\Hone}.
\]
This gives \eqref{eq:local-b-size}.
The modulation equations
\eqref{eq:general-scalar-normal-form}
then yield
\eqref{eq:local-bc-normal}.
Since $b\le C\delta$, choosing $\delta$
sufficiently small gives
\[
c'\ge e_0b-Cb^2
\ge\frac{e_0}{2}b.
\]
Moreover, $|c|\le Cb$ and \(|b'|\le e_0|c|+Cb^2\le Cb\),
proving \eqref{eq:local-monotone-c}.
Finally, for $J=[\alpha,\beta]\subset\tau(I)$,
\[
\begin{aligned}
\int_\alpha^\beta b(\tau)\dd\tau
\le\frac{2}{e_0}
\int_\alpha^\beta c'(\tau)\dd\tau=\frac{2}{e_0}
\bigl(c(\beta)-c(\alpha)\bigr)
\le C\delta,
\end{aligned}
\]
where the last inequality follows from
$|c|\le C\delta$.
This proves \eqref{eq:local-distance-L1}.
\end{proof}

\subsection{Ejection}

We use the hyperbolic ejection mechanism of the radial argument
in \cite[Proposition 3.7]{NR}. The estimate below includes the
translation parameters and places a long exit segment in fixed
NLS coordinates.

\begin{lemma}[Ejection in fixed coordinates]
\label{lem:exit-fixed-entry}
Fix a sufficiently small radius $\delta>0$.
Suppose a strong solution $u$ satisfies $E(u)=E(W)$ and \(\|\nabla u(p)\|_{L^2}<\|\nabla W\|_{L^2}\),
stays in $d_{\mathcal W}(u)\le\delta$ on $[p,q]$, and satisfies
$d_{\mathcal W}(u(q))=\delta$. Let
\[
 r=|\xi_1(p)|=\min_{t\in[p,q]}|\xi_1(t)|,
\]
and assume that $r/\delta$ is sufficiently small. Put
\[
 \lambda_p=\lambda(p),\qquad s=\lambda_p^2(t-q),
 \qquad s_p=\lambda_p^2(p-q),
\]
and freeze the symmetry parameters at $p$:
\[
 \widetilde u(s,y)=e^{-i\theta(p)}\lambda_p^{-\frac{1}{2}}
 u\left(q+\lambda_p^{-2}s,\ a(p)+\lambda_p^{-1}y\right).
\]
Then $\widetilde u=W+h$ solves the same NLS on $[s_p,0]$, and
\begin{equation}\label{eq:exit-fixed-exponential}
 \|h(s)\|_{\Hone}\le C\delta e^{e_0s},\qquad
 \delta e^{e_0s_p}\sim r.
\end{equation}
The constants are independent of $r,p,q$.
\end{lemma}

\begin{proof}
Use $b,c$ from Lemma~\ref{lem:local-integrated-distance}, normalize
the rescaled time by $\tau(q)=0$, and write $\tau_p=\tau(p)$.
Since $b(\tau_p)=r$ is a minimum, $b'(\tau_p)\ge0$ and hence
$c(\tau_p)\ge-Cr^2$. Monotonicity gives $c(\tau)\ge-Cr^2$,
while $b(\tau)\ge r$. Thus, for small $r$,
\[
 \frac{1}{2}b\le\alpha:=b+c\le Cb.
\]
The scalar equations yield
\[
 \alpha'=e_0\alpha+O(b^2),\qquad
 (\log\alpha)'=e_0+O(b).
\]
Integrating and using \eqref{eq:local-distance-L1} gives
\[
 b(\tau)\sim r e^{e_0(\tau-\tau_p)}.
\]
Since $b(0)\sim\delta$, it follows that
\[
 b(\tau)+\|v(\tau)\|_{\Hone}\le C\delta e^{e_0\tau},
 \qquad r\sim\delta e^{e_0\tau_p}.
\]

Integrating the modulation speeds from $\tau_p$ to $\tau$ yields
\[
 |\theta(\tau)-\theta(p)|
 +\left|\log\frac{\lambda(\tau)}{\lambda_p}\right|
 +|\lambda_p(a(\tau)-a(p))|
 \le C\delta e^{e_0\tau}.
\]
Here the scale ratio is first bounded above and below to estimate
$\lambda_p|a_\tau|\le(\frac{\lambda_p}{\lambda})Cb$.
Smoothness of the ground state orbit in $\Hone$, together with
isometry of the symmetry action on $v$, gives
\[
 \|h(s(\tau))\|_{\Hone}\le C\delta e^{e_0\tau}.
\]
Moreover,
\[
 \frac{\mathrm{d}s}{\mathrm{d}\tau}
 =\left(\frac{\lambda_p}{\lambda(\tau)}\right)^2
 =1+O(\delta e^{e_0\tau}).
\]
Since $s=0$ at $\tau=0$, integration backwards from the exit gives
$|s(\tau)-\tau|\le C\delta$. Substitution in the preceding estimates
proves \eqref{eq:exit-fixed-exponential}.
\end{proof}

\subsection{Comparison with the special solutions}
To identify the limiting exit data, we compare these segments with the special solutions, parametrized by their unstable coefficient at time zero. Both kinetic branches enter this comparison; the kinetic inequality will select the lower branch in the final argument.

We denote by $W^+$ the special solution constructed in \cite{DM}
on the upper kinetic branch, with
\(\|\nabla W^+(t)\|_{L^2}>\|\nabla W\|_{L^2}\) and
$W^+(t)\to W$ in $\Hone$ as $t\to\infty$.

\begin{lemma}
\label{lem:ancient-construction}
There exist $a_0>0$ and $C>0$ such that, for every
$a\in(-a_0,a_0)$, there is a radial function
$H_a\in C((-\infty,0];\Hone)$ for which $W+H_a$
solves \eqref{eq:nls} on $(-\infty,0]$, with
$H_0=0$, and
\begin{equation}\label{eq:exit-unstable-curve}
a_+(H_a(0))=a,\quad
\|H_a(t)\|_{\Hone}\le C|a|e^{e_0t},
\quad t\le0.
\end{equation}
The map $a\mapsto H_a(0)$ is continuous in $\Hone$,
and
\[
E(W+H_a(t))=E(W),\quad t\le0.
\]
For $a\ne0$, the solution $W+H_a$ is a time
translate of either $\overline{W^-(-t)}$ or
$\overline{W^+(-t)}$. If
\[
\|\nabla(W+H_a(0))\|_{L^2}
<\|\nabla W\|_{L^2},
\]
the former occurs.
\end{lemma}

\begin{proof}
The construction in \cite[Proposition 6.3, (6.8), and Section 6.3]{DM}
gives the two conjugate time-reversed perturbations $v^1,v^2$ with
\[
 v^j(t)=c_j e^{e_0t}Y_++O_{\Hone}(e^{\frac{3e_0t}{2}})
 \quad(t\to-\infty),\qquad c_1c_2<0.
\]
Put $\alpha_j(t)=a_+(v^j(t))$. The fixed-test spectral equation gives
\[
 \alpha_j'=e_0\alpha_j+O(\|v^j\|_{\Hone}^2)
           =e_0\alpha_j+O(\alpha_j^2).
\]
Thus, sufficiently far into the past, each $\alpha_j$ is strictly
monotone and parametrizes one of the two small intervals adjacent
to zero. For small nonzero $a$, choose the corresponding $j$ and
$T(a)$ with $\alpha_j(T(a))=a$, and set
$H_a(t)=v^j(t+T(a))$ for $t\le0$.
The asymptotic gives $e^{e_0T(a)}\sim|a|$ and hence
\eqref{eq:exit-unstable-curve}. Strong time continuity and the
inverse parametrization give continuity for $a\ne0$; the same
bound gives continuity at zero. Energy and the kinetic branches
are those of the two constructed special solutions.
\end{proof}

The next estimate treats all five zero modes through the center
group. It will be used both at an exit from the tube and to identify an
exponentially convergent threshold solution.

\begin{lemma}[Mixed-boundary comparison]
\label{lem:local-mixed-boundary}
Fix $0<\eta<e_0$. There are $\varepsilon_\eta>0$ and $C_\eta<\infty$
such that the following holds for every $s<0$. Let $h$ and $k$ solve
\eqref{eq:exit-fixed-equation} on $[s,0]$ and satisfy
\[
 \sup_{t\in[s,0]}
 (\|h(t)\|_{\dot H^1}+\|k(t)\|_{\dot H^1})\le\varepsilon_\eta,
 \qquad P_+(h-k)(0)=0.
\]
Then
\begin{equation}\label{eq:formal-mixed-terminal}
 \|h(0)-k(0)\|_{\dot H^1}
 \le C_\eta e^{-\eta s}\|h(s)-k(s)\|_{\dot H^1}.
\end{equation}
\end{lemma}

\begin{proof}
Set $w=h-k$ and $F=\Ncal(h)-\Ncal(k)$. Then
$\partial_t w=\mathcal Aw+F$ and $P_+w(0)=0$.
Partition $[s,0]$ into intervals
\[
s=t_0<t_1<\cdots<t_m=0,\quad J_j=[t_j,t_{j+1}],
\]
with $1/2\le |J_j|\le1$, except when $m=1$.
By \eqref{eq:exit-local-linear} and
\eqref{eq:exit-nonlinear}, the small-data
estimates give
\[
\|h\|_{X(J_j)}+\|k\|_{X(J_j)}
\le C\varepsilon_\eta,
\]
and hence
\begin{equation}\label{eq:mixed-local-forcing}
\|F\|_{N(J_j)}
\le C\varepsilon_\eta\|w\|_{X(J_j)}.
\end{equation}

Decompose \(w=P_+w+P_-w+P_cw\).
By Duhamel's formula and $P_+w(0)=0$,
\[
\begin{aligned}
P_+w(t)
&=-\int_t^0 e^{e_0(t-r)}P_+F(r)\dd r,\\
P_-w(t)
&=e^{-e_0(t-s)}P_-w(s)
+\int_s^t e^{-e_0(t-r)}P_-F(r)\dd r,\\
P_cw(t)
&=e^{(t-s)\mathcal A}P_cw(s)
+\int_s^t e^{(t-r)\mathcal A}P_cF(r)\dd r.
\end{aligned}
\]
Applying \eqref{eq:exit-local-linear}
to the projected equations, forward for
$P_-$ and $P_c$ and backward for $P_+$,
and using \eqref{eq:exit-center-coercivity},
we obtain 
\begin{equation}\label{eq:mixed-slab-estimate}
\begin{aligned}
\|w\|_{X(J_j)}
&\le C(1+t_j-s)\|w(s)\|_{\Hone}+C\sum_{i=0}^j
(1+t_j-t_i)\|F\|_{N(J_i)}\\
&\quad +C\sum_{i=j}^{m-1}
e^{-e_0(t_i-t_j)}\|F\|_{N(J_i)}.
\end{aligned}
\end{equation}
Fix $0<\eta<e_0$. By
\eqref{eq:mixed-local-forcing},
\[
\begin{aligned}
\|F\|_{N(J_i)}
&\le C\varepsilon_\eta e^{\eta(t_i-s)}
\max_{0\le j<m}
\left\{
e^{-\eta(t_j-s)}
\|w\|_{X(J_j)}
\right\}.
\end{aligned}
\]
Since the intervals have lengths comparable
to one,
\[
\begin{aligned}
\sup_j\sum_{i=0}^j
(1+t_j-t_i)e^{-\eta(t_j-t_i)}
\le C_\eta,    \qquad
\sup_j\sum_{i=j}^{m-1}
e^{-(e_0-\eta)(t_i-t_j)}
\le C_\eta.
\end{aligned}
\]
Also,
\[
\sup_{r\ge0}(1+r)e^{-\eta r}<\infty.
\]
Multiplying \eqref{eq:mixed-slab-estimate} by
$e^{-\eta(t_j-s)}$, taking the maximum over $j$,
and choosing $\varepsilon_\eta$ sufficiently small,
we obtain
\[
\max_{0\le j<m}
e^{-\eta(t_j-s)}\|w\|_{X(J_j)}
\le C_\eta\|w(s)\|_{\Hone}.
\]
Since $0\in J_{m-1}$ and $t_{m-1}\le0$,
\[
\begin{aligned}
\|w(0)\|_{\Hone}
\le\|w\|_{X(J_{m-1})}\le C_\eta e^{-\eta s}\|w(s)\|_{\Hone}.
\end{aligned}
\]
This proves \eqref{eq:formal-mixed-terminal}.
\end{proof}

\begin{corollary}[Exponential uniqueness]
\label{cor:local-exponential-uniqueness}
Let $u=W+h$ be a solution of NLS on
$(-\infty,0]$. Suppose that, for some
$C,\beta>0$,
\[
\|h(t)\|_{\dot H^1}
\le Ce^{\beta t},
\quad t\le0.
\]
Then there exist $T<0$ and a sufficiently
small $a\in\mathbb R$ such that
\[
a=a_+(h(T)),
\quad
u(t)=W+H_a(t-T),
\quad t\le T.
\]
\end{corollary}
\begin{proof}
Fix $0<\eta<\min\{\beta,e_0\}$.
For $T<0$, set \(g(t)=h(t+T)\) and \(a=a_+(g(0))\). 
Since $a_+$ is bounded on $\Hone$, it holds that
\[
 |a|\lesssim\|g(0)\|_{\Hone}
 \lesssim e^{\beta T}.
\]
Choose $T<0$ sufficiently negative that
$|a|<a_0$, where $a_0$ is given by
Lemma~\ref{lem:ancient-construction}. Then
\[
\sup_{t\le0}
\bigl(\|g(t)\|_{\Hone}
+\|H_a(t)\|_{\Hone}\bigr)
\lesssim e^{\beta T}.
\]
Decreasing $T$ further if necessary, this
quantity is smaller than $\varepsilon_\eta$. 
Since $P_+(g-H_a)(0)=0$,
Lemma~\ref{lem:local-mixed-boundary} gives
\[
\begin{aligned}
\|g(0)-H_a(0)\|_{\dot H^1}
\le C_\eta e^{-\eta s}
\|g(s)-H_a(s)\|_{\dot H^1}\le C_\eta\bigl(
e^{\beta T}e^{(\beta-\eta)s}
+|a|e^{(e_0-\eta)s}
\bigr).
\end{aligned}
\]
Letting $s\to-\infty$ yields
$g(0)=H_a(0)$.
Uniqueness of the NLS then gives 
\[
h(t)=H_a(t-T),\quad t\le T.
\]
\end{proof}

\section{Proof of the threshold classification}
\label{sec:formal-classification}

We combine threshold compactness with the local analysis to prove
eventual trapping, exponential convergence and identification with $W^-$.

\begin{proposition}[Eventual trapping]
\label{prop:lower-eventual-trapping}
Let $u$ be a forward-global strong solution with
\[
E(u)=E(W),\quad
\|\nabla u(t)\|_{L^2}<\|\nabla W\|_{L^2},
\quad t\ge0.
\]
Suppose that $d_{\Wcal}(u(t_n))\to0$
for some $t_n\to\infty$.
If $u$ does not scatter forward, then,
for every sufficiently small $\delta>0$,
there exists $T_\delta>0$ such that
\[
d_{\Wcal}(u(t))<\delta,\quad t\ge T_\delta.
\]
\end{proposition}
\begin{proof}
Suppose the conclusion fails. Fix a sufficiently
small $\delta>0$, 
for all sufficiently large $n$, let $q_n>t_n$
be the first time such that \(d_{\Wcal}(u(q_n))=\delta\).
Choose $p_n\in[t_n,q_n]$ so that \(r_n=|\xi_1(p_n)|
=\min_{t\in[t_n,q_n]}|\xi_1(t)|
\).
By Lemma~\ref{lem:general-first-exit-coercivity},
\[
0<r_n\le C d_{\Wcal}(u(t_n))\to0.
\]
In particular, $p_n<q_n$ for all sufficiently large $n$.
By Lemma~\ref{lem:exit-fixed-entry}, after fixing
the symmetry parameters at $p_n$, we obtain
solutions $W+h_n$ on $[s_n,0]$ satisfying
\[
\|h_n(s)\|_{\Hone}\le C\delta e^{e_0s},
\qquad
\delta e^{e_0s_n}\sim r_n.
\]
Hence $s_n\to-\infty$.
Set $a_n=a_+(h_n(0))$. Since $|a_n|\le C\delta$,
Lemma~\ref{lem:ancient-construction} gives
an ancient solution $W+H_{a_n}$ with
\[
a_+(H_{a_n}(0))=a_n,
\qquad
\|H_{a_n}(s)\|_{\Hone}\le C\delta e^{e_0s}.
\]
Thus \(P_+(h_n-H_{a_n})(0)=0\)
and
\[
\|h_n(s_n)-H_{a_n}(s_n)\|_{\Hone}
\le C\delta e^{e_0s_n}.
\]
Taking $\delta$ sufficiently small, we may apply
Lemma~\ref{lem:local-mixed-boundary} with
$\eta=e_0/2$ to obtain
\[
\begin{aligned}
\|h_n(0)-H_{a_n}(0)\|_{\Hone}
\le C e^{-\frac{e_0s_n}{2}}
\|h_n(s_n)-H_{a_n}(s_n)\|_{\Hone}\le C\delta e^{\frac{e_0s_n}{2}}
\to0.
\end{aligned}
\]

Passing to a subsequence, assume $a_n\to a_*$.
By the continuity of $a\mapsto H_a(0)$,
the normalized exit data converge strongly
in $\Hone$ to
\[
F_*=W+H_{a_*}(0).
\]
Since the normalization preserves the energy,
kinetic norm, and distance to $\Wcal$,
\[
E(F_*)=E(W),\qquad
\|\nabla F_*\|_{L^2}\le \|\nabla W\|_{L^2},\qquad
d_{\Wcal}(F_*)=\delta.
\]
If $\|\nabla F_*\|_{L^2}= \|\nabla W\|_{L^2}$, energy identity gives
$\|F_*\|_{L^6}=\|W\|_{L^6}$. Equality in the sharp Sobolev
inequality would imply $F_*\in\Wcal$, a contradiction.
Therefore
\[
\|\nabla F_*\|_{L^2} < \|\nabla W\|_{L^2},
\qquad a_*\ne0.
\]

By Lemma~\ref{lem:ancient-construction}, the solution
through $F_*$ is a time translate of
$\overline{W^-(-t)}$, which scatters forward.
Long-time stability then implies that the
normalized solutions with initial data
$W+h_n(0)$ also scatter forward for all
sufficiently large $n$.
Since these are fixed symmetry transforms
of $u$, this contradicts the forward
nonscattering assumption.
\end{proof}

\begin{proposition}[Rigidity of trapped solutions]
\label{prop:lower-trapping-rigidity}
There exists $\delta_*>0$ with the following
property. Let $u$ be a strong solution, global
for $t\ge t_0$, satisfying
\[
 E(u)=E(W),\qquad
 \|\nabla u(t_0)\|_{L^2}
 <\|\nabla W\|_{L^2},
 \qquad
 \sup_{t\ge t_0}d_{\mathcal W}(u(t))
 \le\delta_*.
\]
Then $u$ is $W^-$ modulo the symmetries
of \eqref{eq:nls}.
\end{proposition}

\begin{proof}
By a fixed symmetry and time translation, we may assume
$t_0=0$, $\lambda(0)=1$, and $a(0)=0$.
The strict kinetic inequality is preserved, so $u$ is
nonstationary. We use the modulation parameters from
Lemma~\ref{lem:general-first-exit-coercivity}, with
$\tau(0)=0$, and the functions $b,c$ from
Lemma~\ref{lem:local-integrated-distance}.

By \eqref{eq:local-distance-L1} and
\eqref{eq:formal-mod-speed},
\[
\int_0^{\tau(T)}b(s)\dd s\le C\delta_*,
\quad
|\log\lambda(T)|\le C\delta_*.
\]
For sufficiently small $\delta_*$,
\[
\frac12\le\lambda(t)\le2,\quad t\ge0.
\]
Since $\mathrm d\tau/\mathrm dt=\lambda^2$,
we have $\tau(t)\to\infty$ and
\[
\int_0^\infty b(\tau)\dd\tau<\infty.
\]

By \eqref{eq:local-monotone-c},
\[
|b'|\le Cb,\quad |c|\le Cb,\quad
c'\ge\frac{e_0}{2}b.
\]
The integrability of $b$ and the bound on $b'$
imply $b(\tau)\to0$, hence $c(\tau)\to0$.
Set
\[
D(\tau)=\int_\tau^\infty b(s)\dd s.
\]
Integrating the inequality for $c'$ yields
\[
D(\tau)\le-\frac{2}{e_0}c(\tau)
\le Cb(\tau).
\]
Thus
\[
D'(\tau)=-b(\tau)\le-C^{-1}D(\tau),
\]
which gives $D(\tau)\le Ce^{-\rho\tau}$
for some $\rho>0$. Since $|b'|\le Cb$,
\[
b(\tau)\le C\int_\tau^{\tau+1}b(s)\dd s
\le Ce^{-\rho\tau}.
\]

By \eqref{eq:formal-mod-speed}, the modulation
parameters converge exponentially to
$\theta_\infty$, $\lambda_\infty>0$, and
$a_\infty$. Since $\|v(\tau)\|_{\Hone}\le Cb(\tau)$
and $\tau(t)\ge t/4$, we obtain
\[
\left\|
u(t)-\mathcal T_{\theta_\infty,
\lambda_\infty,a_\infty}W
\right\|_{\Hone}
\le Ce^{-\rho t/4}.
\]

After applying a fixed symmetry, we obtain
a solution $U=W+q$ satisfying
\[
\|q(t)\|_{\Hone}\le Ce^{-\beta t},
\quad t\ge0,
\]
for some $\beta>0$.
Choose $T_1$ sufficiently large and set
\[
h(s)=\overline{U(T_1-s)}-W,\quad s\le0.
\]
Then
\[
\|h(s)\|_{\Hone}
\le Ce^{-\beta T_1}e^{\beta s}.
\]
By Corollary~\ref{cor:local-exponential-uniqueness},
there exist $S<0$ and a sufficiently small
$a\in\mathbb R$ such that
\[
 h(s)=H_a(s-S),\qquad s\le S.
\]
Lemma~\ref{lem:ancient-construction} then shows
that $U$ is either stationary or a time
translate of $W^-$ or $W^+$. The strict kinetic
inequality excludes the stationary solution and
the $W^+$ branch. Thus, by uniqueness, $u$
coincides with $W^-$ modulo symmetries on its
lifespan.
\end{proof}

\begin{proof}[Proof of Theorem~\ref{thm:classification}]
When $\|\nabla u_0\|_{L^2}=\|\nabla W\|_{L^2}$,
the energy identity gives \( \|u_0\|_{L^6}^6
 =\|\nabla u_0\|_{L^2}^2=G\).
Equality in the sharp Sobolev inequality implies
$u_0\in\Wcal$. By uniqueness, $u$ is a stationary
element of $\Wcal$.

Suppose now that \(\|\nabla u_0\|_{L^2}<\|\nabla W\|_{L^2}\).
By Proposition~\ref{prop:lt-global-ap}, $u$ is global.
If $u$ does not scatter forward,
Corollary~\ref{cor:lt-sequential-ground}
gives a sequence $t_n\to\infty$ such that \(d_{\Wcal}(u(t_n))\to0\).
Proposition~\ref{prop:lower-eventual-trapping},
with the trapping radius chosen sufficiently small,
and Proposition~\ref{prop:lower-trapping-rigidity}
then show that $u$ coincides with $W^-$ modulo
symmetries. In particular, $u$ scatters backward.

Suppose instead that $u$ scatters forward.
If it does not scatter backward, applying the
preceding argument to $\overline{u(-t)}$ shows
that $u$ coincides with $\overline{W^-(-t)}$
modulo symmetries. Otherwise, $u$ scatters in
both time directions.
\end{proof}

The classification also gives scattering in both time directions
when $E(u_0)=E(W)$, $\|\nabla u_0\|_{L^2} < \|\nabla W\|_{L^2}$, and either $u_0\in H^1(\R^3)$
or $e^{-i\theta}u_0$ is real-valued for some constant $\theta$.
Indeed, $W^-(t)\notin L^2$ at every finite time in dimension three
\cite[Remark 1.1]{DM}, which excludes the two exceptional orbits
for finite-mass data. In the real-data case, remove the constant
phase. Uniqueness gives $u(-t)=\overline{u(t)}$, so scattering in
one direction, guaranteed by Theorem~\ref{thm:classification},
implies scattering in the other.

\section*{Declarations}
\subsection*{Acknowledgments}
M. Chen and H. Wang were supported by the National Key R\&D Program of China
(2025YFA1018500) and the Beijing Natural Science Foundation (No.~1262019). 

\subsection*{Data availability}
No datasets were generated or analyzed for this article.

\subsection*{Competing interests}
The authors declare no competing interests.

\subsection*{Use of AI}
During the preparation of this manuscript, the authors used ChatGPT (OpenAI) to assist with checking and refining certain mathematical proofs, improving the mathematical exposition, language editing, and LaTeX formatting. The authors assume full responsibility for all content.

\end{document}